\documentclass[11pt]{article}
\usepackage{amsfonts}
\usepackage{amssymb}
\usepackage{color}
\newtheorem{theo}{Theorem}[section]
\newtheorem{lem}[theo]{Lemma}

\newcommand{\mysection}[1]{\section{#1} \setcounter{equation}{0}}
\newcommand{\proof}{{\sc Proof.} \quad}
\newcommand{\proofc}{{\sc Proof} \ }
\newcommand{\be}{\begin{equation} \label}
\newcommand{\ee}{\end{equation}}
\newcommand{\bea}{\begin{eqnarray}\label}
\newcommand{\eea}{\end{eqnarray}}
\newcommand{\bas}{\begin{eqnarray*}}
\newcommand{\eas}{\end{eqnarray*}}
\newcommand{\bit}{\begin{itemize}}
\newcommand{\eit}{\end{itemize}}
\newcommand{\qed}{\hfill$\Box$ \vskip.2cm}
\newcommand{\nn}{\nonumber}
\newcommand{\R}{\mathbb{R}}

\newcommand{\pO}{\partial\Omega}

\newcommand{\eps}{\varepsilon}

\newcommand{\abs}{\\[5pt]}
\newcommand{\Abs}{\\[5mm]}

\newcommand{\cred}{\color{red}}
\newcommand{\tm}{T_{max}}
\newcommand{\io}{\int_\Omega}
\newcommand{\bom}{\overline{\Omega}}
\newcommand{\ouz}{\overline{u}_0}
\newcommand{\ovz}{\overline{v}_0}
\newcommand{\oy}{\overline{y}}
\newcommand{\ws}{w_\star}
\newcommand{\ow}{\overline{w}}
\newcommand{\F}{{\mathcal{F}}}
\newcommand{\D}{{\mathcal{D}}}
\begin{document}
\title{Analysis of a one-dimensional forager-exploiter model}
\author
{
Youshan Tao\footnote{taoys@dhu.edu.cn}\\
{\small Department of Applied Mathematics, Dong Hua University,}\\
{\small Shanghai 200051, P.R.~China}
 \and
Michael Winkler\footnote{michael.winkler@math.uni-paderborn.de}\\
{\small Institut f\"ur Mathematik, Universit\"at Paderborn,}\\
{\small 33098 Paderborn, Germany} }
\date{}
\maketitle
\begin{abstract}
\noindent
  The system
  \bas
  \left\{ \begin{array}{l}
    u_t= u_{xx} - \chi_1 (uw_x)_x, \\[1mm]
    v_t = v_{xx} - \chi_2 (vu_x)_x,  \\[1mm]
    w_t = dw_{xx} - \lambda (u+v)w - \mu w + r,
    \end{array} \right.
  \eas
  is considered in a bounded real interval, with positive parameters $\chi_1,\chi_2,d,\lambda$ and $\mu$,
  and with $r \ge 0$.
  Proposed to describe social interactions within mixed forager-exploiter groups, this model
  extends classical one-species chemotaxis-consumption systems by additionally accounting for a second
  taxis mechanism coupled to the first in a consecutive manner. \abs
  It is firstly shown that for all suitably regular initial data $(u_0, v_0, w_0)$, an
  associated Neumann-type initial-boundary value problem
  possesses a globally defined bounded classical solution.
  Moreover, it is asserted that this solution stabilizes to a spatially homogeneous equilibrium at an exponential rate
  under a smallness condition on $\min\{\io u_0, \io v_0\}$ that appears to be consistent with predictions
  obtained from formal stability analysis.\abs
\end{abstract}
\newpage
\mysection{Introduction}
It is well-known that social interactions in mixed-species groups may lead to
rich spatial patterns, already in some cases of populations consisting of merely two fractions,
and the ambition to understand fundamental principles for such complex dynamics
has been attracting considerable interest during the past decade in the literature near the borderline regions
between, inter alia, biology, behavioral sciences
and mathematics (\cite{giraldeau_caraco}, \cite{tania_pnas}, \cite{guttal}, \cite{bellomo_soler_M3AS2012}).
With the macroscopic formation of shearwater flocks through attraction to kittiwake foragers in Alaska
forming a paradigmatic example, a rather clearly traceable type of social interplay, in the comparatively simple case
of only two involved species,       
is constituted by so-called forager-exploiter interaction:
The members of a first population, the ``foragers", search for food by directly  
moving upward gradients of the nutrient concentration; contrary to
this, as ``exploiters" the individuals of a second population pursue
a more indirect strategy by rather orienting their movement toward
regions of higher forager population densities.\abs
As a macroscopic model for the spatio-temporal evolution of the population densities $u=u(x,t)$ and $v=v(x,t)$ of foragers
and scroungers in such contexts,
additionally accounting for the nutrient density $w=w(x,t)$ as a third unknown
the authors in \cite{tania_pnas} propose the parabolic PDE system
\be{00}
    \left\{ \begin{array}{l}
    u_t =\Delta u - \chi_1 \nabla\cdot (u\nabla w), \\[1mm]
    v_t = \Delta v -\chi_2 \nabla \cdot (v\nabla u), \\[1mm]
    w_t = d\Delta v - \lambda (u+v)w - \mu w + r,
    \end{array} \right.
\ee
with positive parameters $\chi_1, \chi_2, d$ and $\lambda$, and with $\mu\ge 0$ and $r\ge 0$.
Besides presupposing random diffusion of all quantities, this model assumes the food resources to be degraded by
both foragers and scroungers upon contact, and to possibly be supplied by an external source and spontaneously decaying.
In its most characteristic part, however, (\ref{00}) accounts for the group-specific strategies of directed motion
by postulating taxis-type cross-diffusion mechanisms to be responsible in both cases.\abs
As a particular     
feature thereby generated, (\ref{00}) contains a sequential coupling of two taxis processes,
which may be expected to considerably increase the mathematical complexity of (\ref{00}) when compared, for instance,
to the
corresponding one-species chemotaxis-consumption system, depending on the application context sometimes also referred
to as prey-taxis system (\cite{kareiva_odell}, \cite{lee_hillen_lewis}, \cite{junping_shi_JDE}),
that is obtained on letting $v\equiv 0$ in (\ref{00}), and that is hence given by
\be{01}
    \left\{ \begin{array}{l}
    u_t =\Delta u - \chi_1 \nabla\cdot (u\nabla w), \\[1mm]
    w_t = d\Delta v - \lambda uw - \mu w + r.
    \end{array} \right.
\ee
In the prototypical case when $\mu=r=0$, namely, the structure of the latter is artless enough so as to allow for a meaningful
energy structure that can be used as a technical basis for a comprehensive theory of global classical well-posedness
in low-dimensional boundary value problems in which the spatial
dimension $n$ satisfies $n\le 2$, of global weak solvability when
$n=3$, and even of asymptotic stabilization toward spatially
homogeneous equilibria whenever $n\le 3$
(\cite{taowin_consumption}). Adaptations of such approaches have
been utilized to address several variants of (\ref{01}), partially
even involving further components and interaction mechanisms
(\cite{wang_zhian}, \cite{junping_shi_JDE}, \cite{win_JDE2017},
\cite{win_ARMA}), but capturing the additional intricacy induced by
the second taxis interaction in (\ref{00}) seems beyond the
abilities of such methods.\abs
In particular, to the best of our knowledge it is yet left completely open by the analytical literature
how far the coupling of the nutrient taxis mechanism from (\ref{01}) to
a further cross-diffusion process sensitive to the gradient of the first population may lead to substantial
destabilization of the tendency toward homogeneity, as known to occur in (\ref{01});
only under some restrictive assumptions on $r$ and the initial data $w|_{t=0}$, in essence ensuring that $w$
remains below a suitably small threshold throughout evolution, it has recently been possible to achieve
some results on global existence of certain generalized
solutions, as well as on their large time stabilization toward constants when moreover $r=r(x,t)$ decays suitably in time
(\cite{win_scrounge}).
In light of the strong singularity-supporting potential of chemotactic cross-diffusion, well-known as a striking
feature e.g.~of the classical Keller-Segel system (\cite{herrero_velazquez}, \cite{win_JMPA}),
already answers to questions from basic theory of global solvability in (\ref{00}) thus seem far from obvious,
even in the simplest case in which $n=1$; in fact, the numerical experiments reported in \cite{tania_pnas} indicate quite a rich
dynamical potential of (\ref{00}) already in such one-dimensional frameworks.\Abs
{\bf Main results.} \quad
The purpose of the present work now consists in making sure that the evident mathematical challenges notwithstanding,
and especially despite the apparent lack of any favorable energy structure, at least the one-dimensional version of
(\ref{00}) does not only allow for a rather comprehensive theory of classical solvability, but beyond this is even accessible
to an essentially exhaustive qualitative analysis in parameter constellations for which formal considerations predict
asymptotic homogenization:
In particular, we shall develop an analytical approach which firstly enables us to
assert global existence of bounded classical solutions for widely arbitrary initial
data, and which secondly is subtle enough so as to allow for a conclusion on large time stabilization toward constant
steady states under an additional smallness assumption on the total population sizes of foragers and scroungers
that quite precisely seems to match a corresponding condition formally identified as essentially necessary and sufficient
therefor by means of a linear stability analysis in \cite{tania_pnas}.\abs
To take this more precise, in a bounded open interval $\Omega\subset
\R$ let us consider the initial-boundary value problem for
(\ref{00}) given by \be{0}
        \left\{ \begin{array}{lcll}
    u_t &=& u_{xx} - \chi_1 (uw_x)_x,
    \qquad & x\in\Omega, \ t>0, \\[1mm]
    v_t &=& v_{xx} - \chi_2 (vu_x)_x,
    \qquad & x\in\Omega, \ t>0, \\[1mm]
    w_t &=& dw_{xx} - \lambda (u+v)w - \mu w + r,
    \qquad & x\in\Omega, \ t>0, \\[1mm]
    & & \hspace*{-14mm}
    u_x=v_x=w_x=0,
    \qquad & x\in\pO, \ t>0, \\[1mm]
    & & \hspace*{-14mm}
    u(x,0)=u_0(x), \quad v(x,0)=v_0(x), \quad w(x,0)=w_0(x),
    \qquad & x\in\Omega,
    \end{array} \right.
\ee
where $\chi_1,\chi_2,d,\lambda$ and $\mu$ are positive constants and $r$ is nonnegative,
and where the initial data are such that
\be{init}
    \left\{ \begin{array}{l}
    u_0\in W^{1,\infty}(\Omega)
    \mbox{ is nonnegative with } u_0\not\equiv 0, \\[1mm]
    v_0\in W^{1,\infty}(\Omega)
    \mbox{ is nonnegative with } v_0\not\equiv 0,
    \qquad \mbox{and that} \\[1mm]
    w_0\in W^{1,\infty}(\Omega)
    \mbox{ is positive in } \bom.
    \end{array} \right.
\ee
Then the first of our main results asserts global existence of bounded classical solutions in the following flavor.
\begin{theo}\label{theo34}
  Let $\Omega\subset\R$ be a bounded open interval, and let
  $\chi_1,\chi_2,d,\lambda$ and $\mu$ be positive and $r$ be nonnegative.
  Then for any choice of $(u_0,v_0,w_0)$ fulfilling (\ref{init}), the problem (\ref{0}) possesses
  a global classical solution $(u,v,w)$ which is uniquely determined by the properties that
  \be{34.01}
    \mbox{$u,v$ and $w$ belong to }
    C^0([0,\infty);W^{1,2}(\Omega)) \cap C^{2,1}(\bom\times (0,\infty)),
  \ee
  and which is such that $u>0$ and $v>0$ in $\bom\times (0,\infty)$ as well as $w>0$ in $\bom\times [0,\infty)$.
  Moreover, this solution is bounded in the sense that there exists $C>0$ such that
  \be{34.1}
    \|u(\cdot,t)\|_{W^{1,2}(\Omega)}
    + \|v(\cdot,t)\|_{W^{1,2}(\Omega)}
    + \|w(\cdot,t)\|_{W^{1,2}(\Omega)}
    \le C
    \qquad \mbox{for all } t>0.
  \ee
\end{theo}
Next intending to identify circumstances under which the diffusion processes in (\ref{0}) are sufficiently strong
so as to warrant relaxation into homogeneous states,
we note that in light of corresponding results on asymptotically diffusive behavior in chemotaxis systems with
essentially superlinear nonlinear ingredients (\cite{cao_optspace}, \cite{cao_lankeit}, \cite{win_JDE})
it seems far from audacious to conjecture
that such dynamics can always be observed if solutions remain small {\em in all their components},
and {\em with respect to suitably fine topologies}, throughout evolution.\abs
In order to approach a more subtle picture in this regard, we recall
that a linear stability analysis detailed in \cite{tania_pnas}
suggests to expect, in the normalized case when $\Omega=(0,1)$ and
the spatial averages $\ouz:=\frac{1}{|\Omega|}\io u_0$ and
$\ovz:=\frac{1}{|\Omega|} v_0$ satisfy $\ouz+\ovz=1$, the relation
\be{cond}
    \frac{8(\lambda+\mu)^2 (d+1)}{\lambda r \chi_1\ouz\ovz}
    +\frac{2(d+1)}{\ovz}  \gtrsim \chi_2
\ee as the decisive condition for prevalence of homogeneity. In
particular, this condition is entirely independent of $w_0$, and
moreover it involves $u_0$ and $v_0$ exclusively through their $L^1$
norms as their biologically best interpretable derivate. If, more
generally, $\ouz+\ovz$ is supposed to remain below a given number in
an arbitrary bounded $\Omega$, then assuming (\ref{cond}) evidently
becomes equivalent to imposing a certain smallness hypothesis on
$\min\{\ouz,\ovz\}$. The second of our main results now provides a
rigorous mathematical counterpart of this formal consideration.
\begin{theo}\label{theo40}
  Suppose that $\Omega\subset\R$ is a bounded open interval, that $\chi_1,\chi_2,d,\lambda$ and $\mu$ are positive, and that
  $r\ge 0$.
  Then for all $M>0$ one can find $\eps(M)>0$ with the property that whenever $u_0, v_0$ and $w_0$ are such that besides
  (\ref{init}) we have
  \be{M}
    \io u_0 + \io v_0 \le M
  \ee
  as well as
  \be{del1}
    \io u_0 \le \eps(M)
    \qquad \mbox{or} \qquad
    \io v_0 \le \eps(M),
  \ee
  there exist $C=C(u_0,v_0,w_0)>0$ and $\alpha=\alpha(u_0,v_0,w_0)>0$ such that the solution $(u,v,w)$ of (\ref{0}) satisfies
  \be{40.1}
    \|u(\cdot,t)-\ouz\|_{L^\infty(\Omega)}
    + \|v(\cdot,t)-\ovz\|_{L^\infty(\Omega)}
    + \|w(\cdot,t)-\ws\|_{L^\infty(\Omega)}
    \le C e^{-\alpha t}
    \qquad \mbox{for all } t>0,
  \ee
  where $\ws$ is the nonnegative constant given by
  \be{ws}
    \ws:=\frac{r}{\lambda(\ouz+\ovz)+\mu}.
  \ee
\end{theo}
\mysection{Local existence and an explicit $L^\infty$ bound for $w$}
The following basic statement on local existence and extensibility can be obtained from standard theory on evolution
systems of parabolic type.
\begin{lem}\label{lem_loc}
  Suppose that $\Omega\subset\R$ is a bounded open interval, that
  $\chi_1,\chi_2,d,\lambda$ and $\mu$ are positive and $r\ge 0$, and that (\ref{init}) holds.
  Then there exist $\tm \in (0,\infty]$ and nonnegative functions $u,v,w$, uniquely determined by the requirement that
  \be{reg}
    \mbox{$u,v$ and $w$ are elements of }
    C^0([0,\tm);W^{1,2}(\Omega)) \cap C^{2,1}(\bom\times (0,\tm)),
  \ee
  which solve (\ref{0}) in the classical sense in $\Omega\times (0,\tm)$, and which are such that
  \be{ext}
    \mbox{if $\tm<\infty$, then} \quad
    \limsup_{t\nearrow\tm} \Big\{ \|u(\cdot,t)\|_{W^{1,q}(\Omega)} + \|v(\cdot,t)\|_{W^{1,q}(\Omega)}
    + \|w(\cdot,t)\|_{W^{1,q}(\Omega)} \Big\}
    = \infty
    \qquad \mbox{for all } q>1.
  \ee
  Moreover, $u>0$ and $v>0$ in $\bom\times (0,\tm)$ and $w>0$ in $\bom\times [0,\tm)$, and we have
  \be{mass}
    \io u(\cdot,t) = \io u_0
    \quad \mbox{and} \quad
    \io v(\cdot,t) = \io v_0
    \qquad \mbox{for all } t\in (0,\tm).
  \ee
\end{lem}
\proof
  Since all eigenvalues of the diffusion matrix
  \bas
    \mathcal {A}(u, v, w) :=\Bigg(
    \begin{array}{ccc}
    1 & 0 & -\chi_1 u\\
    -\chi_2 v & 1 & 0\\
    0 & 0 &d
    \end{array} \Bigg)
  \eas
  are positive, the general theory on local existence and maximal extension from \cite{amann}
  is applicable to (\ref{0}); in particular, the existence of a uniquely
  determined maximal classical solution follows from \cite[Theorems
  14.4 and 14.6]{amann}, and the extensibility criterion (\ref{ext})
  is ensured by \cite[Theorem 15.5]{amann}.
  Moreover, the strong maximum principle along with our assumptions $u_0\not\equiv 0$ and $v_0\not\equiv
  0$ in (\ref{init}) yield positivity of $u$ and $v$ in $\bom\times(0,\tm)$,
  whereas a simple comparison argument shows that $w(x, t)\ge \min_{x\in\bom}w_0(x) e^{-\mu t}>0$ in $\bom\times [0,\tm)$
  thanks to our assumption that $w_0>0$.
  Finally, the mass conservation properties in (\ref{mass}) immediately result from integration in
  the first equation and the second equation in (\ref{0}), respectively.
\qed
Constituting another basic but important feature of (\ref{0}), the following pointwise bound on $w$ is an immediate
consequence of the maximum principle.
\begin{lem}\label{lem26}
  We have
  \be{26.1}
    \|w(\cdot,t)\|_{L^\infty(\Omega)}
    \le \frac{r}{\mu} + \|w_0\|_{L^\infty(\Omega)} e^{-\mu t}
    \qquad \mbox{for all } t\in (0,\tm).
  \ee
\end{lem}
\proof
  We let $\ow(x,t):= y(t)$ for $x\in\bom$ and $t\ge 0$, where
  \bas
    y(t):=y_0 e^{-\mu t} + \frac{r}{\mu} (1-e^{-\mu t}),
    \qquad t\ge 0,
  \eas
  denotes the solution of $y'(t) + \mu y(t)=r$, $t>0$, with $y(0)=y_0:=\|w_0\|_{L^\infty(\Omega)}$.
  Then $\ow(\cdot,0) = y_0 \ge w(\cdot,0)$ in $\Omega$ as well as $\frac{\partial \ow}{\partial\nu}=0$ on
  $\pO\times (0,\infty)$.
  As moreover, by nonnegativity of $\lambda, u$ and $v$,
  \bas
    \ow_t - d\ow_{xx} + \lambda(u+v) w + \mu \ow - r
    &\ge& \ow_t - d\ow_{xx} + \mu \ow - r \\
    &\ge& y' + \mu y -r=0
    \qquad \mbox{in } \Omega\times (0,\tm),
  \eas
  by means of a comparison argument we conclude that $w\le \ow$ in $\Omega\times (0,\tm)$, and that thus
  \bas
    w(x,t) \le y_0 e^{-\mu t} + \frac{r}{\mu}
    \qquad \mbox{for all $x\in\Omega$ and } t\in (0,\tm),
  \eas
  which is equivalent to (\ref{26.1}).
\qed
\mysection{Further estimates. Linking regularity to the size of $\io u_0+\io v_0$}\label{sect3}
The goal of this section is to reveal further regularity properties of the above local solution
which on the one hand will allow for its global extension, but which on the other will also prepare our
subsequent qualitative analysis.
For this purpose, the dependence of the obtained estimates on the averages $\ouz$ and $\ovz$ as well as on
their sum, as appearing in (\ref{M}), will carefully be traced throughout this section.
\subsection{A bound for $w_x$ in $L^q$}
In a first step we shall employ parabolic smoothing estimates to see that the mere mass conservation properties
in (\ref{mass}) entail $L^q$ bounds for the chemotactic gradient acting in the first equation from (\ref{0}).
In not requiring any restriction other than that $q$ be finite, we here
already make essential use of our assumption on the spatial framework to be one-dimensional.
\begin{lem}\label{lem25}
  There exists $\alpha>0$ such that for all $M>0$ and any $q>1$ one can find $K(M,q)>0$
  with the property that whenever $u_0,v_0$ and $w_0$ satisfy (\ref{init}) as well as (\ref{M}), there exists
  $C=C(u_0,v_0,w_0)>0$ such that
  \be{25.1}
    \|w_x(\cdot,t)\|_{L^q(\Omega)}
    \le K(M,q) + C e^{-\alpha t}
    \qquad \mbox{for all } t\in (0,\tm).
  \ee
\end{lem}
\proof
  Relying on known regularization features of the Neumann heat semigroup $(e^{\sigma\Delta})_{\sigma\ge 0}$ on $\Omega$
  with $\Delta:=(\cdot)_{xx}$ (\cite{win_JDE}), let us fix $c_1(q)>0$ and $c_2(q)>0$ such that for all $t>0$,
  \be{25.2}
    \|\partial_x e^{dt\Delta} \varphi\|_{L^q(\Omega)}
    \le c_1(q)\|\varphi\|_{W^{1,\infty}(\Omega)}
    \qquad \mbox{for all } \varphi\in W^{1,\infty}(\Omega)
  \ee
  and
  \be{25.3}
    \|\partial_x e^{dt\Delta} \varphi\|_{L^q(\Omega)}
    \le c_2(q) (1+t^{-1+\frac{1}{2q}}) \|\varphi\|_{L^1(\Omega)}
    \qquad \mbox{for all } \varphi\in C^0(\bom).
  \ee
  Then representing $w$ according to
  \bas
    w(\cdot,t)
    = e^{t(d\Delta-\mu)} w_0
    + \int_0^t e^{(t-s)(d\Delta-\mu)} \Big\{ -\lambda (u(\cdot,s)+v(\cdot,s)) w(\cdot,s) + r \Big\} ds,
    \qquad t\in (0,\tm),
  \eas
  we can combine (\ref{25.2}) with (\ref{25.3}) to estimate
  \bea{25.4}
     & & \hspace*{-16mm} \|w_x(\cdot,t)\|_{L^q(\Omega)}\nn\\
    &\le& c_1(q) e^{-\mu t} \|w_0\|_{W^{1,\infty}(\Omega)}\nn\\
    & &   + c_2(q) \int_0^t e^{-\mu(t-s)} \Big(1+(t-s)^{-1+\frac{1}{2q}}\Big)
        \Big\| -\lambda(u(\cdot,s)+v(\cdot,s)) w(\cdot,s) + r \Big\|_{L^1(\Omega)} ds \nn\\
    &\le& c_1(q) e^{-\mu t} \|w_0\|_{W^{1,\infty}(\Omega)} \nn\\
    & & + c_2(q) \lambda \int_0^t e^{-\mu(t-s)}  \Big(1+(t-s)^{-1+\frac{1}{2q}}\Big) \cdot
        \Big\{ \|u(\cdot,s)\|_{L^1(\Omega)} + \|v(\cdot,s)\|_{L^1(\Omega)} \Big\} \|w(\cdot,s)\|_{L^\infty(\Omega)} ds
        \nn\\
    & & + c_2(q) r |\Omega| \int_0^t e^{-\mu(t-s)} \Big(1+(t-s)^{-1+\frac{1}{2q}}\Big) ds
    \qquad \mbox{for all } t\in (0,\tm).
  \eea
  Here using (\ref{mass}) along with (\ref{M}) and Lemma \ref{lem26}, we find that
  \bas
    & & \hspace*{-20mm}
    c_2(q) \lambda \int_0^t e^{-\mu(t-s)}  \Big(1+(t-s)^{-1+\frac{1}{2q}}\Big) \cdot
    \Big\{ \|u(\cdot,s)\|_{L^1(\Omega)} + \|v(\cdot,s)\|_{L^1(\Omega)} \Big\} \|w(\cdot,s)\|_{L^\infty(\Omega)} ds \\
    &\le& c_2(q)\lambda \bigg\{ \io u_0 + \io v_0\bigg\} \cdot
    \int_0^t e^{-\mu(t-s)}  \Big(1+(t-s)^{-1+\frac{1}{2q}}\Big) \cdot \Big\{ \frac{r}{\mu} + \|w_0\|_{L^\infty(\Omega)} e^{-\mu s} \Big\} ds \\
    & \le & c_2(q) \lambda \bigg\{ \io u_0+ \io v_0\bigg\} \cdot
    \bigg\{ \frac{r}{\mu} \cdot \Big(\frac{2}{\mu} +2q\Big) + \|w_0\|_{L^\infty(\Omega)} \cdot
    (te^{-\mu t} +2q t^\frac{1}{2q}e^{-\mu t})\bigg\} \\
    &\le& \frac { 2c_2(q)\lambda M r(1+q \mu)}{\mu^2}
    + c_2(q) \lambda M \|w_0\|_{L^\infty(\Omega)} \cdot  (te^{-\mu t} +2q t^\frac{1}{2q}e^{-\mu t})
    \qquad \mbox{for all } t\in (0,\tm)
  \eas
  because $ \int_0^t e^{-\mu(t-s)} \Big(1+(t-s)^{-1+\frac{1}{2q}}\Big) ds=\int_0^t e^{-\mu s}(1+s^{-1+\frac{1}{2q}})ds
  \le \int_0^\infty e^{-\mu s}(1+s^{-1+\frac{1}{2q}})ds \le \frac{2}{\mu}+2 q $ and
  $ \int_0^t e^{-\mu t} \Big(1+(t-s)^{-1+\frac{1}{2q}}\Big) ds =te^{-\mu t} +2q t^\frac{1}{2q}e^{-\mu t}$, whence moreover estimating
  $te^{-\mu t} \le \frac{2}{\mu} e^{-\frac{\mu}{2}t}$
  as well as $2q t^\frac{1}{2q}e^{-\mu t} \le \frac{2}{\mu} e^{-\frac{\mu}{2}t}$ for $t>0$,
  from (\ref{25.4}) we infer that
  \bas
    \|w_x(\cdot,t)\|_{L^q(\Omega)}
    &\le& \frac{ 2c_2(q)\lambda M r(1+q\mu)}{\mu^2}
    + \frac{2c_2(q) r(1+q\mu) |\Omega|}{\mu} \\
    & & + \bigg\{ c_1(q) \|w_0\|_{W^{1,\infty}(\Omega)}
    + \frac{4c_2(q)\lambda M \|w_0\|_{L^\infty(\Omega)}}{\mu} \bigg\} \cdot e^{-\frac{\mu}{2} t}
    \qquad \mbox{for all } t\in (0,\tm),
  \eas
  which is precisely of the claimed form with suitably chosen $K(M,q)>0$ and $C=C(u_0,v_0,w_0)>0$, and with
  $\alpha:=\frac{\mu}{2}$.
\qed
\subsection{Estimating $u$ in $L^\infty$ and $u_x$ in a spatio-temporal $L^2$ norm}
Through a standard testing procedure performed using the first equation in (\ref{0}), when applied to $q:=2$
the latter has a first consequence on regularity of $u$ as well as its gradient.
\begin{lem}\label{lem23}
  There exists $\alpha>0$ with the property that given any $M>0$ one can choose $K(M)>0$
  such that if (\ref{init}) and (\ref{M}) hold, there exists
  $C=C(u_0,v_0,w_0)>0$ such that
  \be{23.1}
    \io u^2(\cdot,t)
    \le K(M) \ouz^2 + C e^{-\alpha t}
    \qquad \mbox{for all } t\in (0,\tm)
  \ee
  and
  \be{23.2}
    \int_t^{t+\tau} \io u_x^2
    \le K(M) \ouz^2 + C e^{-\alpha t}
    \qquad \mbox{for all } t\in (0,\tm-\tau),
  \ee
  where $\tau:=\min\{1,\frac{1}{2}\tm\}$.
\end{lem}
\proof
  Let us first apply Lemma \ref{lem25} to $q:=2$ to fix positive constants $\alpha$ and $k_1(M)$ such that whenever
  (\ref{init}) and (\ref{M}) hold, one can find $c_1=c_1(u_0,v_0,w_0)>0$ such that
  \be{23.22}
    \|w_x(\cdot,t)\|_{L^2(\Omega)} \le k_1(M) + c_1 e^{-\alpha t}
    \qquad \mbox{for all } t\in (0,\tm),
  \ee
  where without loss of generality we may assume that $k_1(M)\ge 1$ and $\alpha<\frac{1}{6}$.
  Apart from that, by means of the Gagliardo-Nirenberg inequality and Young's inequality we can choose $c_2>0$, $c_3>0$ and
  $c_4>0$ such that
  \be{23.3}
    \|\varphi\|_{L^\infty(\Omega)}
    \le c_2 \|\varphi_x\|_{L^2(\Omega)}^\frac{2}{3} \|\varphi\|_{L^1(\Omega)}^\frac{1}{3}
    + c_2\|\varphi\|_{L^1(\Omega)}
    \qquad \mbox{for all } \varphi\in W^{1,2}(\Omega),
  \ee
  that
  \be{23.33}
    \frac{1}{2} \|\varphi\|_{L^2(\Omega)}^2
    \le \frac{1}{4} \|\varphi_x\|_{L^2(\Omega)}^2
    + c_3 \|\varphi\|_{L^1(\Omega)}^2
    \qquad \mbox{for all } \varphi\in W^{1,2}(\Omega),
  \ee
  and that
  \be{23.4}
    c_2 \chi_1  |\Omega|^\frac{1}{3} ab \le \frac{1}{8} a^\frac{6}{5} + c_4 b^6
    \qquad \mbox{for all $a\ge 0$ and } b\ge 0.
  \ee
  Now assuming (\ref{init}) and (\ref{M}), we test the first equation in (\ref{0}) by $u$ to see using the Cauchy-Schwarz
  inequality as well as (\ref{23.22}), (\ref{23.3}), (\ref{23.33}) and (\ref{mass}), that
  \bea{23.5}
    & & \hspace*{-14mm}
    \frac{1}{2} \frac{d}{dt} \io u^2 + \frac{1}{2} \io u^2 + \io u_x^2 \nn\\
    &=& \chi_1 \io uu_x w_x
    + \frac{1}{2} \io u^2 \nn\\
    &\le& \chi_1 \|u\|_{L^\infty(\Omega)} \|u_x\|_{L^2(\Omega)} \|w_x\|_{L^2(\Omega)}
    + \frac{1}{2} \io u^2 \nn\\
    &\le& \chi_1 \|u\|_{L^\infty(\Omega)} \|u_x\|_{L^2(\Omega)} \cdot
    \Big\{ k_1(M) + c_1 e^{-\alpha t} \Big\}
    + \frac{1}{2} \io u^2 \nn\\
    &\le& c_2 \chi_1 \|u_x\|_{L^2(\Omega)}^\frac{5}{3} \|u\|_{L^1(\Omega)}^\frac{1}{3} \cdot
    \Big\{ k_1(M) + c_1 e^{-\alpha t} \Big\}
    + c_2 \chi_1 \|u_x\|_{L^2(\Omega)} \|u\|_{L^1(\Omega)} \cdot
    \Big\{ k_1(M) + c_1 e^{-\alpha t} \Big\} \nn\\
    & & + \frac{1}{4} \|u_x\|_{L^2(\Omega)}^2
    + c_3\|u\|_{L^1(\Omega)}^2 \nn\\
    &=& c_2 \chi_1 |\Omega|^\frac{1}{3} \ouz^\frac{1}{3} \|u_x\|_{L^2(\Omega)}^\frac{5}{3} \cdot
    \Big\{ k_1(M) + c_1 e^{-\alpha t} \Big\}
    + c_2 \chi_1 |\Omega| \ouz  \|u_x\|_{L^2(\Omega)}\cdot
    \Big\{ k_1(M) + c_1 e^{-\alpha t} \Big\} \nn\\
    & & + \frac{1}{4} \|u_x\|_{L^2(\Omega)}^2
    + c_3 |\Omega|^2 \ouz^2
    \qquad \mbox{for all } t\in (0,\tm),
  \eea
  where by (\ref{23.4}),
  \bas
    c_2 \chi_1 |\Omega|^\frac{1}{3} \ouz^\frac{1}{3} \|u_x\|_{L^2(\Omega)}^\frac{5}{3} \cdot
    \Big\{ k_1(M) + c_1 e^{-\alpha t} \Big\}
    \le \frac{1}{8} \|u_x\|_{L^2(\Omega)}^2
    + c_4 \ouz^2 \cdot
    \Big\{ k_1(M) + c_1 e^{-\alpha t} \Big\}^6
  \eas
 and where once more by Young's inequality,
  \bas
    c_2 \chi_1 |\Omega| \ouz \|u_x\|_{L^2(\Omega)} \cdot
    \Big\{ k_1(M) + c_1 e^{-\alpha t} \Big\}
    \le \frac{1}{8} \|u_x\|_{L^2(\Omega)}^2
    + 2c_2^2 \chi_1^2 |\Omega|^2 \ouz^2 \cdot
    \Big\{ k_1(M) + c_1 e^{-\alpha t} \Big\}^2
  \eas
  for all $t\in (0,\tm)$.
  Since $k_1(M)\ge 1$ and thus
  \bas
    \Big\{ k_1(M) + c_1 e^{-\alpha t} \Big\}^2
    \le \Big\{ k_1(M) + c_1 e^{-\alpha t} \Big\}^6
     \le
    32 \cdot \Big\{ k_1^6(M) + c_1^6 e^{-6\alpha t} \Big\}
    \qquad \mbox{for all } t>0,
  \eas
  from (\ref{23.5}) we altogether obtain that
  \be{23.6}
    \frac{d}{dt} \io u^2 + \io u^2 + \io u_x^2
    \le k_2(M) \ouz^2
    + c_5 e^{-6\alpha t}
    \qquad \mbox{for all } t\in (0,\tm)
  \ee
  with
  \bas
    k_2(M):=2c_3 |\Omega|^2 + 64 c_4 k_1^6(M) + 128 c_2^2 \chi_1^2 |\Omega|^2 k_1^6(M)
  \eas
  and
  \bas
    c_5\equiv c_5(u_0,v_0,w_0):=64  c_1^6 c_4 \ouz^2 + 128 c_1^6 c_2^2 \chi_1^2 |\Omega|^2 \ouz^2.
  \eas
  Using that $6\alpha<1$, we may invoke Lemma \ref{lem21} to firstly conclude from (\ref{23.6}) that
  \be{23.7}
    \io u^2(\cdot,t)
    \le \bigg\{ \io u_0^2 + \frac{c_5}{1-6\alpha} \bigg\} \cdot e^{-6\alpha t}
    + k_2(M) \ouz^2
    \qquad \mbox{for all } t\in (0,\tm),
  \ee
  and that thus (\ref{23.1}) holds with evident choices of the constants therein.
  After that, by direct integration of (\ref{23.6}) we see that since $\tau\le 1$,
  \bas
    \int_t^{t+\tau} \io u_x^2
    &\le& \io u^2(\cdot,t)
    + k_2(M) \ouz^2
    + c_5 \int_t^{t+\tau} e^{-6\alpha s} ds \\
    &\le& \io u^2(\cdot,t)
    + k_2(M) \ouz^2
    + c_5 e^{-6\alpha t}
    \qquad \mbox{for all } t\in (0,\tm-\tau),
  \eas
  which shows that (\ref{23.7}) also entails (\ref{23.2}).
\qed
By again going back to Lemma \ref{lem25}, in light of the outcome of Lemma \ref{lem23} we can now once more employ
heat semigroup estimates to actually improve the
topological setting in (\ref{23.1}) so as to involve the respective $L^\infty$ norm.
\begin{lem}\label{lem24}
  There exists $\alpha>0$ such that if $M>0$, then one can fix $K(M)>0$
  such that under the assumptions (\ref{init}) and (\ref{M}) it is possible to find
  $C=C(u_0,v_0,w_0)>0$ such that
  \be{24.1}
    \|u(\cdot,t)\|_{L^\infty(\Omega)}
    \le K(M) \ouz + C e^{-\alpha t}
    \qquad \mbox{for all } t\in (0,\tm).
  \ee
\end{lem}
\proof
  We begin by employing Lemma \ref{lem23} and Lemma \ref{lem25} to take $\alpha_1\in (0,1)$ and $\alpha_2\in (0,1)$
  such that given $M>0$ we can find $k_1(M)>0$ and $k_2(M)>0$ with the property that if (\ref{init}) and (\ref{M}) hold,
  then with some $c_i=c_i(u_0,v_0,w_0)>0$, $i\in\{1,2\}$, we have
  \bas
    \|u(\cdot,t)\|_{L^2(\Omega)}
    \le k_1(M) \ouz + c_1 e^{-\alpha_1 t}
    \qquad \mbox{for all } t\in (0,\tm)
  \eas
  and
  \bas
    \|w_x(\cdot,t)\|_{L^4(\Omega)}
    \le k_2(M) + c_2 e^{-\alpha_2 t}
    \qquad \mbox{for all } t\in (0,\tm)
  \eas
  and hence, by the H\"older inequality,
  \bea{24.2}
    \|u(\cdot,t) w_x(\cdot,t)\|_{L^\frac{4}{3}(\Omega)}
    &\le& \|u(\cdot,t)\|_{L^2(\Omega)} \|w_x(\cdot,t)\|_{L^4(\Omega)} \nn\\
    &\le& k_1(M) k_2(M) \ouz
    + c_1 k_2(M) e^{-\alpha_1 t}
    + c_2 k_1(M) \ouz e^{-\alpha_2 t}
    + c_1 c_2 e^{ -(\alpha_1+\alpha_2)t} \nn\\
    &\le& k_1(M) k_2(M) \ouz
    + c_3 e^{-\alpha t}
    \qquad \mbox{for all } t\in (0,\tm)
  \eea
  with $\alpha:=\min\{\alpha_1,\alpha_2\}$ and $c_3=c_3(u_0,v_0,w_0):=c_1 k_2(M) + c_2 k_1(M) \ouz + c_1 c_2$.
  Next, parabolic smoothing estimates (\cite{win_JDE},  \cite{FIWY}) provide $c_4>0$ and $c_5>0$ such that for all $t>0$,
  \bas
    \|e^{t\Delta} \varphi_x\|_{L^\infty(\Omega)}
    \le c_4 (1+t^{-\frac{7}{8}}) \|\varphi\|_{L^\frac{4}{3}(\Omega)}
    \qquad \mbox{for all $\varphi\in C^1(\bom)$ such that $\varphi_x=0$ on } \pO
  \eas
  and
  \bas
    \|e^{t\Delta} \varphi\|_{L^\infty(\Omega)}
    \le c_5 (1+t^{-\frac{1}{2}}) \|\varphi\|_{L^1(\Omega)}
    \qquad \mbox{for all } \varphi\in C^0(\bom),
  \eas
  so that henceforth assuming (\ref{init}) and (\ref{M}) for some $M>0$, and
  rewriting the first equation in (\ref{0}) in the form
  \bas
    u_t -u_{xx}+u = -\chi_1 (uw_x)_x + u
    \qquad \mbox{in } \Omega\times (0,\tm),
  \eas
  by means of an associated variation-of-constants representation we may estimate
  \bea{24.3}
    \|u(\cdot,t)\|_{L^\infty(\Omega)}
    &=& \bigg\| e^{t(\Delta-1)} u_0
    - \chi_1 \int_0^t e^{(t-s)(\Delta-1)} \partial_x \Big(u(\cdot,s) w_x(\cdot,s)\Big) ds
    + \int_0^t e^{(t-s)(\Delta-1)} u(\cdot,s) ds \bigg\|_{L^\infty(\Omega)} \nn\\
    &\le& e^{-t} \|e^{t\Delta} u_0\|_{L^\infty(\Omega)}
    + c_4 \chi_1 \int_0^t e^{-(t-s)} \Big(1+(t-s)^{-\frac{7}{8}}\Big) \|u(\cdot,s) w_x(\cdot,s)\|_{L^\frac{4}{3}(\Omega)} ds
        \nn\\
    & & + c_5 \int_0^t e^{-(t-s)} \Big(1+(t-s)^{-\frac{1}{2}}\Big) \|u(\cdot,s)\|_{L^1(\Omega)} ds
    \qquad \mbox{for all } t\in (0,\tm).
  \eea
  Here by the maximum principle and the fact that $\alpha\le 1$,
  \be{24.4}
    e^{-t} \|e^{t\Delta} u_0\|_{L^\infty(\Omega)}
    \le \|u_0\|_{L^\infty(\Omega)} e^{-t}
    \le \|u_0\|_{L^\infty(\Omega)} e^{ -\alpha t}
    \qquad \mbox{for all } t>0,
  \ee
  whereas (\ref{mass}) ensures that
  \bea{24.5}
    c_5 \int_0^t e^{-(t-s)} \Big(1+(t-s)^{-\frac{1}{2}}\Big) \|u(\cdot,s)\|_{L^1(\Omega)} ds
    &=& c_5 |\Omega| \ouz \int_0^t e^{-(t-s)} \Big( 1+(t-s)^{-\frac{1}{2}}\Big) ds \nn\\
    &\le& c_5 c_6 |\Omega| \ouz
    \qquad \mbox{for all } t\in (0,\tm)
  \eea
  with $c_6:=\int_0^\infty e^{-\sigma} (1+\sigma^{-\frac{1}{2}}) d\sigma<\infty$.
  Moreover, thanks to (\ref{24.2}) we have
  \bea{24.6}
    & & \hspace*{-20mm}
    c_4 \chi_1 \int_0^t e^{-(t-s)} \Big(1+(t-s)^{-\frac{7}{8}}\Big) \|u(\cdot,s) w_x(\cdot,s)\|_{L^\frac{4}{3}(\Omega)} ds
        \nn\\
    &\le& c_4 \chi_1 k_1(M) k_2(M) \ouz \int_0^t e^{-(t-s)} \Big(1+(t-s)^{-\frac{7}{8}}\Big) ds \nn\\
    & & + c_3 c_4 \chi_1 \int_0^t e^{-(t-s)} \Big(1+(t-s)^{-\frac{7}{8}}\Big) e^{-\alpha s} ds
    \qquad \mbox{for all } t\in (0,\tm),
  \eea
  where
  \bas
    \int_0^t e^{-(t-s)} \Big(1+(t-s)^{-\frac{7}{8}}\Big) ds
    \le c_7:=\int_0^\infty e^{-\sigma} (1+\sigma^{-\frac{7}{8}}) d\sigma
    \qquad \mbox{for all } t>0,
  \eas
  and where
  \bas
    \int_0^t e^{-(t-s)} \Big(1+(t-s)^{-\frac{7}{8}}\Big) e^{-\alpha s} ds
    = e^{-\alpha t} \int_0^t e^{-(1-\alpha)\sigma} (1+\sigma^{-\frac{7}{8}}) d\sigma
    \le c_8 e^{-\alpha t}
    \qquad \mbox{for all } t>0
  \eas
  with $c_8:=\int_0^\infty e^{-(1-\alpha)\sigma} (1+\sigma^{-\frac{7}{8}}) d\sigma$ being finite thanks to our restriction
  that $\alpha<1$.\abs
  Inserting (\ref{24.4})-(\ref{24.6}) into (\ref{24.3}) thus shows that for all $t\in (0,\tm)$,
  \bas
    \|u(\cdot,t)\|_{L^\infty(\Omega)}
    \le \Big\{ c_5 c_6 |\Omega| + c_4 c_7 \chi_1 k_1(M) k_2(M)\Big\} \cdot \ouz
    + \Big\{ \|u_0\|_{L^\infty(\Omega)} + c_3 c_4 c_8 \chi_1 \Big\} \cdot e^{-\alpha t},
  \eas
  and therefore establishes (\ref{24.1}) upon the observation that $c_4, c_5, c_6$ and $c_7$ do not depend on our particular
  choice of $u_0, v_0$ and $w_0$.
\qed
\subsection{Space-time $L^2$ bounds for $v$ and for $w_{xx}$}
Now unlike in the analysis of (\ref{01}), for globally extending our solution the bounds obtained in Lemma \ref{lem24}
and Lemma \ref{lem25} seem yet insufficient: In view of the second equation in (\ref{0})
it seems that for the detection of appropriate estimates for the second solution component,
further information on the respectively relevant cross-diffusive gradient $u_x$ seems in order.
To prepare our derivation thereof in the next section, let us here provide some preliminary bounds on $v$, $v_x$ and $w_{xx}$
useful for that purpose.\abs
We begin with a basic space-time integrability feature of $(\ln (v+1))_x$ which by another straightforward testing procedure
can be seen to be quite a direct consequence of
our present knowledge on $u_x$ from Lemma \ref{lem23}.
\begin{lem}\label{lem27}
  One can find $\alpha>0$ with the property that to any $M>0$ there corresponds some $K(M)>0$
  such that if (\ref{init}) and (\ref{M}) hold, then there exists
  $C=C(u_0,v_0,w_0)>0$ fulfilling
  \be{27.1}
    \int_t^{t+\tau} \io \frac{v_x^2}{(v+1)^2}
    \le K(M) + C e^{-\alpha t}
    \qquad \mbox{for all } t\in (0,\tm-\tau),
  \ee
  where $\tau:=\min\{1,\frac{1}{2}\tm\}$.
\end{lem}
\proof
  We multiply the second equation in (\ref{0}) by $\frac{1}{v+1}$ and integrate by parts to see that due to Young' inequality,
  \bas
    \frac{d}{dt} \io \ln (v+1)
    &=& \io \frac{v_x^2}{(v+1)^2}
    - \chi_2 \io \frac{v}{(v+1)^2} u_x v_x \\
    &\ge& \frac{1}{2} \io \frac{v_x^2}{(v+1)^2}
    - \frac{\chi_2^2}{2} \io \frac{v^2}{(v+1)^2} u_x^2 \\
    &\ge& \frac{1}{2} \io \frac{v_x^2}{(v+1)^2}
    - \frac{\chi_2^2}{2} \io u_x^2
    \qquad \mbox{for all } t\in (0,\tm).
  \eas
  As $0\le \ln (\xi+1) \le \xi$ for all $\xi\ge 0$, further integration shows that thanks to (\ref{mass}),
  \bas
    \frac{1}{2} \int_t^{t+\tau} \io \frac{v_x^2}{(v+1)^2}
    &\le& \io \ln \Big(v(\cdot,t+\tau)+1\Big) - \io \ln \Big(v(\cdot,t)+1\Big)
    + \frac{\chi^2}{2} \int_t^{t+\tau} \io u_x^2 \\
    &\le& \io v_0
    + \frac{\chi^2}{2} \int_t^{t+\tau} \io u_x^2
    \qquad \mbox{for all  }t\in (0,\tm,\tau),
  \eas
  so that (\ref{27.1}) becomes a consequence of Lemma \ref{lem23}.
\qed
Again thanks to the fact that the considered setting is one-dimensional, a simple interpolation argument shows that
the above entails a space time bound on $v$ itself, rather than on the quantity $\ln (v+1)$ addressed in Lemma \ref{lem27}.
\begin{lem}\label{lem28}
  There exists $\alpha>0$ such that whenever $M>0$, one can pick $K(M)>0$
  such that if (\ref{init}) and (\ref{M}) hold, then with some
  $C=C(u_0,v_0,w_0)>0$ we have
  \be{28.1}
    \int_t^{t+\tau} \io v^2
    \le K(M) + C e^{-\alpha t}
    \qquad \mbox{for all } t\in (0,\tm-\tau),
  \ee
  where again $\tau:=\min\{1,\frac{1}{2}\tm\}$.
\end{lem}
\proof
  According to the one-dimensional Gagliardo-Nirenberg inequality, we can fix $c_1>0$ such that
  \bas
    \|\varphi\|_{L^4(\Omega)}^4
    \le c_1 \|\varphi_x\|_{L^1(\Omega)}^2 \|\varphi\|_{L^2(\Omega)}^2
    + c_1\|\varphi\|_{L^2(\Omega)}^4
    \qquad \mbox{for all } \varphi\in W^{1,1}(\Omega),
  \eas
  which when applied to $\sqrt{v(\cdot,t)+1}$, $t\in (0,\tm)$, shows that since
  $\|\sqrt{v+1}\|_{L^2(\Omega)}^2=\io v_0+|\Omega|$ for all $t\in (0,\tm)$ by (\ref{mass}),
  \bas
    \io v^2
    \le \io (v+1)^2
    &=& \|\sqrt{v+1}\|_{L^4(\Omega)}^4 \\
    &\le& c_1\|\partial_x \sqrt{v+1}\|_{L^1(\Omega)}^2 \|\sqrt{v+1}\|_{L^2(\Omega)}^2
    + c_1\|\sqrt{v+1}\|_{L^2(\Omega)}^4 \\
    &=& \frac{c_1}{4} \cdot \bigg\{ \io v_0+|\Omega| \bigg\} \cdot \bigg\{ \io \frac{|v_x|}{\sqrt{v+1}} \bigg\}^2
    + c_1 \cdot \bigg\{ \io v_0+|\Omega| \bigg\}^2
  \eas
  for all $t\in (0,\tm)$.
  Once more in view of (\ref{mass}), using the Cauchy-Schwarz inequality we see that herein
  \bas
    \bigg\{ \io \frac{|v_x|}{\sqrt{v+1}} \bigg\}^2
    &\le& \bigg\{ \io (v+1) \bigg\} \cdot \io \frac{v_x^2}{(v+1)^2} \\
    &\le& \bigg\{ \io v_0 + |\Omega| \bigg\} \cdot \io \frac{v_x^2}{(v+1)^2}
    \qquad \mbox{for all } t\in (0,\tm),
  \eas
  whence altogether, after a time integration,
  \bas
    \int_t^{t+\tau} \io v^2
    \le c_1 \cdot \bigg\{ \io v_0 + |\Omega| \bigg\}^2 \cdot \bigg\{ \frac{1}{4} \io \frac{v_x^2}{(v+1)^2} +1 \bigg\}
    \qquad \mbox{for all } t\in (0,\tm-\tau)
  \eas
  due to the fact that $\tau\le 1$.
  The claimed statement thus readily results from Lemma \ref{lem27}.
\qed
Thus having at hand spatio-temporal integral estimates for $v$ and, through e.g.~Lemma \ref{lem24}, also for $u$,
we have collected sufficient regularity information on all the source terms in the third equation from (\ref{0}),
when considered as a semilinear heat equation, so as to obtain the following as the outcome of a further standard testing
process.
\begin{lem}\label{lem29}
  One can find $\alpha>0$ in such a way that for each $M>0$ there exists $K(M)>0$
  such that assuming (\ref{init}) and (\ref{M}) entails that with some
  $C=C(u_0,v_0,w_0)>0$,
  \be{29.1}
    \int_t^{t+\tau} \io w_{xx}^2
    \le \frac{K(M)}{\tau} + C e^{-\alpha t}
    \qquad \mbox{for all } t\in (0,\tm-\tau),
  \ee
  where once more $\tau:=\min\{1,\frac{1}{2}\tm\}$.
\end{lem}
\proof
  By means of Lemma \ref{lem23} and Lemma \ref{lem28}, we can find $\alpha_1>0$ and $\alpha_2>0$ such that given any $M>0$
  one can pick $k_1(M)>0$ and $k_2(M)>0$ such that whenever (\ref{init}) and (\ref{M}) hold, for the corresponding
  solution of (\ref{0}) we have
  \be{29.2}
    \io u^2 \le k_1(M) + c_1 e^{-\alpha_1 t}
    \qquad \mbox{for all } t\in (0,\tm)
  \ee
  and
  \be{29.3}
    \int_t^{t+\tau} \io v^2 \le k_2(M) + c_2 e^{-\alpha_2 t}
    \qquad \mbox{for all } t\in (0,\tm-\tau)
  \ee
  with some $c_i=c_i(u_0,v_0,w_0)>0$, $i\in\{1,2\}$, and $\tau=\min\{ 1,\frac{1}{2}\tm\}$.\abs
  Now indeed assuming (\ref{init}) and (\ref{M}), we use $w_{xx}$ as a test function for the third equation in (\ref{0})
  to see that by Young's inequality,
  \bas
    \frac{1}{2} \frac{d}{dt} \io w_x^2 + \mu \io w_x^2 + d \io w_{xx}^2
    &=& \lambda \io uww_{xx} + \lambda \io vww_{xx} \\
    &\le& \frac{d}{2} \io w_{xx}^2
    + \frac{\lambda^2}{d} \io u^2 w^2
    + \frac{\lambda^2}{d} \io v^2 w^2 \\
    &\le& \frac{d}{2} \io w_{xx}^2
    + \frac{\lambda^2}{d} \|w\|_{L^\infty(\Omega)}^2 \cdot \bigg\{ \io u^2 + \io v^2 \bigg\}
  \eas
  for all $t\in (0,\tm)$ and hence
  \be{29.4}
    \frac{d}{dt} \io w_x^2 + 2\mu \io w_x^2 + d \io w_{xx}^2
    \le \frac{2\lambda^2}{d} \|w\|_{L^\infty(\Omega)}^2 \cdot \bigg\{ \io u^2 + \io v^2 \bigg\}
    \qquad \mbox{for all } t\in (0,\tm).
  \ee
  Here aiming at an application of Lemma \ref{lem211}, we pick any $\alpha>0$ such that $\alpha<2\mu$ and
  $\alpha\le \min\{\alpha_1,\alpha_2\}$, and combine (\ref{29.2}) and (\ref{29.3}) with the outcome of Lemma \ref{lem26}
  and Young's inequality to estimate
  \bas
    & & \hspace*{-20mm}
    \int_t^{t+\tau} \Bigg\{ \frac{2\lambda^2}{d} \|w(\cdot,s)\|_{L^\infty(\Omega)}^2 \cdot
    \bigg\{ \io u^2(\cdot,s) + \io v^2(\cdot,s) \bigg\} \Bigg\} ds \\
    &\le& \frac{2\lambda^2}{d} \int_t^{t+\tau} \bigg\{ \frac{r}{\mu} + \|w_0\|_{L^\infty(\Omega)} e^{-\mu t} \bigg\}^2
    \cdot \bigg\{ \io u^2(\cdot,s)+\io v^2(\cdot,s) \bigg\} ds \\
    &\le& \frac{4\lambda^2}{d} \cdot \bigg\{ \frac{r^2}{\mu^2} + \|w_0\|_{L^\infty(\Omega)}^2 e^{-2\mu t} \bigg\}
    \cdot \bigg\{ \int_t^{t+\tau} \io u^2 + \int_t^{t+\tau} \io v^2 \bigg\} \\
    &\le& \frac{4\lambda^2}{d} \cdot \bigg\{ \frac{r^2}{\mu^2} + \|w_0\|_{L^\infty(\Omega)}^2 e^{-2\mu t} \bigg\}
    \cdot \Big\{ k_1(M)+k_2(M) + c_1 e^{-\alpha_1 t} + c_2 e^{-\alpha_2 t} \Big\}
  \eas
  for all $t\in (0,\tm-\tau)$,
  where we have used that $\tau\le 1$ and that hence $\int_t^{t+\tau} e^{-\beta s} ds \le e^{-\beta t}$ for all $t>0$
  and any $\beta>0$. Since $\alpha \le \min\{\alpha_1,\alpha_2,2\mu\}$, this readily implies that
  \be{29.6}
    \int_t^{t+\tau} \Bigg\{ \frac{2\lambda^2}{d} \|w(\cdot,s)\|_{L^\infty(\Omega)}^2 \cdot
    \bigg\{ \io u^2(\cdot,s) + \io v^2(\cdot,s) \bigg\} \Bigg\} ds
    \le k_3(M) + c_4 e^{-\alpha t}
    \qquad \mbox{for all } t\in (0,\tm-\tau)
  \ee
  with $k_3(M):=\frac{4\lambda^2 r^2}{d {\cred \mu^2}} (k_1(M)+k_2(M))$ and
  \bas
    c_4\equiv c_4(u_0,v_0,w_0)
    :=\frac{4\lambda^2 r^2}{d\mu^2} (c_1+c_2)
    + \frac{4\lambda^2}{d} \|w_0\|_{L^\infty(\Omega)}^2 (k_1(M) + k_2(M) + c_1 + c_2).
  \eas
  Upon employing Lemma \ref{lem211}, we thus obtain that (\ref{29.4}) firstly entails the inequality
  \bas
    \io w_x^2
    \le \frac{k_3(M)}{2\mu \tau} + k_3(M) + c_5 e^{-\alpha t}
    \qquad \mbox{for all } t\in (0,\tm)
  \eas
  if we let
  \bas
    c_5\equiv c_5(u_0,v_0,w_0)
    := \frac{e^\alpha}{\tau} \cdot \bigg\{ \io w_{0x}^2 +  c_4 + \frac{c_4}{2\mu-\alpha} + k_3(M) \bigg\}
    + c_4 e^\alpha.
  \eas
  Thereupon, directly integrating (\ref{29.4}) shows that again due to (\ref{29.6}),
  \bas
    d\int_t^{t+\tau} \io w_{xx}^2
    &\le& \io w_x^2(\cdot,t)
    + k_3(M) + c_4 e^{-\alpha t} \\
    &\le& \frac{k_3(M)}{2\mu \tau} + 2k_3(M) + (c_4+c_5) e^{-\alpha t}
    \qquad \mbox{for all } t\in (0,\tm-\tau),
  \eas
  which yields (\ref{29.1}) upon again recalling that $\tau\le 1$.
\qed
\subsection{Estimating $u_x$ in $L^2$}
Thanks to Lemma \ref{lem29}, we now have appropriate information on the coefficient functions
$a(x,t):=-\chi_1 w_x$ and $b(x,t):=-\chi_1 w_{xx}$ in the identity $u_t=u_{xx}+a(x,t)u_x + b(x,t)u$ to see that again due
to a variational argument, $u_x$ indeed enjoys the following integrability features which go substantially beyond
those obtained in Lemma \ref{lem23}.
\begin{lem}\label{lem22}
  There exists $\alpha>0$ such that for arbitrary $M>0$ it is possible to choose $K(M)>0$ in such a way that
  whenever (\ref{init}) and (\ref{M}) are satisfied, there exists
  $C=C(u_0,v_0,w_0)>0$ such that writing $\tau:=\min\{1,\frac{1}{2}\tm\}$ we have
  \be{22.1}
    \io u_x^2(\cdot,t)
    \le \frac{K(M)}{\tau^2} + C e^{-\alpha t}
    \qquad \mbox{for all } t\in (0,\tm)
  \ee
  and
  \be{22.2}
    \int_t^{t+\tau} \io u_{xx}^2
    \le \frac{K(M)}{\tau^2} + C e^{-\alpha t}
    \qquad \mbox{for all } t\in (0,\tm-\tau).
  \ee
\end{lem}
\proof
  As a consequence of Lemma \ref{lem24} and Lemma \ref{lem29}, we may pick $\alpha_1>0$ and $\alpha_2>0$ such that for all $M>0$
  we can find $k_1(M)>0$ and $k_2(M)>0$ with the property that under the assumptions (\ref{init}) and (\ref{M})  one may fix
  $c_1=c_1(u_0,v_0,w_0)>0$ and $c_2=c_2(u_0,v_0,w_0)>0$ fulfilling
  \be{22.3}
    \|u(\cdot,t)\|_{L^\infty(\Omega)}^2
    \le k_1(M) + c_1 e^{-\alpha_1 t}
    \qquad \mbox{for all } t\in (0,\tm)
  \ee
  and
  \be{22.4}
    \int_t^{t+\tau} \io w_{xx}^2 \le \frac{k_2(M)}{\tau} + c_2 e^{-\alpha_2 t}
    \qquad \mbox{for all } t\in (0,\tm-\tau).
  \ee
  Apart from that, we combine the Gagliardo-Nirenberg inequality with Young's inequality to obtain $c_3>0$ and $c_4>0$ such that
  \be{22.5}
    \|\varphi_x\|_{L^4(\Omega)}^2 \le c_3 \|\varphi_{xx}\|_{L^2(\Omega)} \|\varphi\|_{L^\infty(\Omega)}
    \qquad \mbox{for all } \varphi \in W^{2,2}(\Omega)
  \ee
  and that
  \be{22.6}
    \io \varphi_x^2 \le \frac{1}{2} \io \varphi_{xx}^2 + c_4 \bigg\{ \io |\varphi| \bigg\}^2
    \qquad \mbox{for all } \varphi \in W^{2,2}(\Omega)
  \ee
  Now assuming (\ref{init}) and (\ref{M}) to be valid for some $M>0$, we integrate by parts in the first equation from (\ref{0})
  and use the Cauchy-Schwarz inequality along with (\ref{22.5}), Young's inequality and (\ref{22.6}) to see that
  for all $t\in (0,\tm)$,
  \bas
    \frac{d}{dt} \io u_x^2 + \io u_x^2 + 2\io u_{xx}^2
    &=& 2\chi_1 \io u_x w_x u_{xx}
    + 2\chi_1 \io uu_{xx} w_{xx}
    + \io u_x^2 \\
    &=& - \chi_1 \io u_x^2 w_{xx}
    + 2\chi_1 \io uu_{xx} w_{xx}
    + \io u_x^2 \\
    &\le& \chi_1 \|u_x\|_{L^4(\Omega)}^2 \|w_{xx}\|_{L^2(\Omega)}
    + 2\chi_1 \|u\|_{L^\infty(\Omega)} \|u_{xx}\|_{L^2(\Omega)} \|w_{xx}\|_{L^2(\Omega)}
    + \io u_x^2 \\
    &\le& (c_3+2) \chi_1\|u\|_{L^\infty(\Omega)} \|u_{xx}\|_{L^2(\Omega)} \|w_{xx}\|_{L^2(\Omega)}
    + \io u_x^2 \\
    &\le& \io u_{xx}^2
    + \frac{(c_3+2)^2 \chi_1^2}{2} \|u\|_{L^\infty(\Omega)}^2 \|w_{xx}\|_{L^2(\Omega)}^2
    + c_4 \bigg\{ \io u_0\bigg\}^2
  \eas
  because of (\ref{mass}).
  In view of the hypothesis (\ref{M}), this shows that abbreviating $c_5:=\frac{(c_3+2)^2 \chi_1^2}{2}$ we have
  \be{22.7}
    \frac{d}{dt} \io  u_x^2 + \io u_x^2 + \io u_{xx}^2
    \le c_5 \|u\|_{L^\infty(\Omega)}^2 \io w_{xx}^2
    + c_4  M^2
    \qquad \mbox{for all } t\in (0,\tm),
  \ee
  where thanks to (\ref{22.3}) and (\ref{22.4}), fixing any $\alpha\in (0,1)$ such that $\alpha\le \min\{\alpha_1,\alpha_2\}$
  we can estimate
  \bea{22.8}
    & & \hspace*{-20mm}
    \int_t^{t+\tau} \bigg\{ c_5\|u(\cdot,s)\|_{L^\infty(\Omega)}^2 \io w_{xx}^2(\cdot,s) + c_4 M^2 \bigg\} ds \nn\\
    &\le& c_5 \cdot \Big\{ k_1(M)+c_1 e^{-\alpha_1 t} \Big\} \cdot
    \Big\{ \frac{k_2(M)}{\tau} + c_2 e^{-\alpha_2 t} \Big\} + c_4 M^2 \nn\\
    &\le& \frac{k_3(M)}{\tau} + c_6 e^{-\alpha t}
    \qquad \mbox{for all } t\in (0,\tm-\tau)
  \eea
  with $k_3(M):=c_5  k_1(M) k_2(M)+ c_4 M^2$ and
  $c_6\equiv c_6(u_0,v_0,w_0):=c_2 c_5 k_1(M) + \frac{c_1 c_5  k_2(M)}{\tau} + c_1 c_2 c_5$.
  As a consequence of Lemma \ref{lem211} and other restriction that $\alpha<1$, from (\ref{22.7}) we thus infer that
  writing $c_7\equiv c_7(u_0,v_0,w_0):=\frac{e^\alpha}{\tau} \cdot \Big\{ \io u_{0x}^2
    +c_6
  + \frac{c_6}{1-\alpha} + \frac{k_3(M)}{\tau}\Big\} + c_6 e^\alpha$ we have
  \be{22.9}
    \io u_x^2
    \le \frac{k_3(M)}{\tau^2} + \frac{k_3(M)}{\tau} + c_7 e^{-\alpha t}
    \qquad \mbox{for all } t\in (0,\tm),
  \ee
  and that hence, by integration of (\ref{22.7}) and again using (\ref{22.8}),
  \bea{22.10}
    \int_t^{t+\tau} \io u_{xx}^2
    &\le& \io u_x^2(\cdot,t) + \frac{k_3(M)}{\tau} + c_6 e^{-\alpha t} \nn\\
    &\le& \frac{k_3(M)}{\tau^2}
    + \frac{2k_3(M)}{\tau}
    + (c_6+c_7) e^{-\alpha t}
    \qquad \mbox{for all } t\in (0,\tm-\tau).
  \eea
  Since $\tau\le 1$, the claimed properties directly result from (\ref{22.9}) and (\ref{22.10}).
\qed
\subsection{An $L^2$ bound for $v_x$}
We can thereby gradually improve our knowledge on the second solution component, firstly addressing $v$ itself
in the course of a further testing procedure:
\begin{lem}\label{lem31}
  There exists $\alpha>0$ such that for all $M>0$ one can fix $K(M)>0$ having the property that
  whenever (\ref{init}) and (\ref{M}) hold, with some
  $C=C(u_0,v_0,w_0)>0$ we have
  \be{31.1}
    \io v^4(\cdot,t)
    \le \frac{K(M) \ovz^4}{\tau^{10}} + C e^{-\alpha t}
    \qquad \mbox{for all } t\in (0,\tm).
  \ee
\end{lem}
\proof
  On the basis of Lemma \ref{lem22}, it is possible to pick $\alpha_1\in (0,1)$ in such a way that given $M>0$ we can choose
  $k_1(M)>0$ which is such that if (\ref{init}) and (\ref{M}) hold,
  \be{31.2}
    \bigg\{ \io u_x^2 \bigg\}^5 \le \frac{k_1(M)}{\tau^{10}} + c_1 e^{-\alpha_1 t}
    \qquad \mbox{for all } t\in (0,\tm)
  \ee
  with some $c_1=c_1(u_0,v_0,w_0)>0$.
  Once more relying on the Gagliardo-Nirenberg inequality and Young's inequality, we furthermore fix $c_2>0$, $c_3>0$ and
  $c_4>0$ such that
  \be{31.3}
    \|\varphi\|_{L^\infty(\Omega)}
    \le c_2 \|\varphi_x\|_{L^2(\Omega)}^\frac{4}{5} \|\varphi\|_{L^\frac{1}{2}(\Omega)}^\frac{1}{5}
    + c_2\|\varphi\|_{L^\frac{1}{2}(\Omega)}
    \qquad \mbox{for all } \varphi\in W^{1,2}(\Omega),
  \ee
  that
  \be{31.4}
    \io \varphi^2 \le \|\varphi_x\|_{L^2(\Omega)}^2 + c_3 \|\varphi\|_{L^\frac{1}{2}(\Omega)}^2
    \qquad \mbox{for all } \varphi\in W^{1,2}(\Omega),
  \ee
  and that
  \be{31.5}
    6c_2\chi_2 ab
    \le a^\frac{10}{9} + c_4 b^{10}
    \qquad \mbox{for all $a\ge 0$ and } b\ge 0.
  \ee
  Then supposing that (\ref{init}) and (\ref{M}) to be satisfied, we use the second equation in (\ref{0}) to see that due to
  the Cauchy-Schwarz inequality and (\ref{mass}), applications of (\ref{31.3}), (\ref{31.5}), (\ref{31.4}) and
  Young's inequality show that for all $t\in (0,\tm)$,
  \bea{31.6}
    \hspace*{-10mm}
    \frac{d}{dt} \io v^4 + \io v^4 + 3\io (v^2)_x^2
    &=& 6\chi_2 \io v^2 u_x (v^2)_x
    + \io v^4 \nn\\
    &\le& 6\chi_2 \|v^2\|_{L^\infty(\Omega)} \|u_x\|_{L^2(\Omega)} \|(v^2)_x\|_{L^2(\Omega)}
    + \io v^4 \nn\\
    &\le& 6c_2\chi_2 \|(v^2)_x\|_{L^2(\Omega)}^\frac{9}{5} \|v^2\|_{L^\frac{1}{2}(\Omega)}^\frac{1}{5} \|u_x\|_{L^2(\Omega)}
        \nn\\
    & & + 6c_2\chi_2 \|(v^2)_x\|_{L^2(\Omega)} \|v^2\|_{L^\frac{1}{2}(\Omega)} \|u_x\|_{L^2(\Omega)}  \nn\\
    & & + \io v^4 \nn\\
    &\le& \|(v^2)_x\|_{L^2(\Omega)}^2
    + c_4 \|v^2\|_{L^\frac{1}{2}(\Omega)}^2 \|u_x\|_{L^2(\Omega)}^{10} \nn\\
    & & + \|(v^2)_x\|_{L^2(\Omega)}^2
    + 9c_2^2 \chi_2^2 \|v^2\|_{L^\frac{1}{2}(\Omega)}^2 \|u_x\|_{L^2(\Omega)}^2 \nn\\
    & & + \|(v^2)_x\|_{L^2(\Omega)}^2
    + c_3 \|v^2\|_{L^\frac{1}{2}(\Omega)}^2 \nn\\
    &=& 3\io (v^2)_x^2
    + \Big\{ c_4 \|u_x\|_{L^2(\Omega)}^{10} + 9c_2^2 \chi_2^2 \|u_x\|_{L^2(\Omega)}^2 + c_3 \Big\}
        \cdot \|v_0\|_{L^1(\Omega)}^4.
  \eea
  Since here, by using Young's inequality and relying on the fact that $\tau\le 1$, we can estimate
  \bas
    & & \hspace*{-20mm}
    \Big\{ c_4 \|u_x\|_{L^2(\Omega)}^{10} + 9c_2^2 \chi_2^2 \|u_x\|_{L^2(\Omega)}^2 + c_3 \Big\}
     \cdot \|v_0\|_{L^1(\Omega)}^4  \\
    &\le& (c_4 + 9c_2^2 \chi_2^2) \|v_0\|_{L^1(\Omega)}^4 \cdot \bigg\{ \io u_x^2 \bigg\}^5
    + (9c_2^2 \chi_2^2 + c_3) \|v_0\|_{L^1(\Omega)}^4 \\
    &\le& \frac{k_2(M) \ovz^4}{\tau^{10}} + c_5 e^{-\alpha_1 t}
    \qquad \mbox{for all } t\in (0,\tm)
  \eas
  with $k_2(M):=\Big\{ (c_4+9c_2^2\chi_2^2) k_1(M) + 9c_2^2 \chi_2^2 + c_3 \Big\} \cdot |\Omega|^4$ and
  $c_5\equiv c_5(u_0,v_0,w_0):=c_1 \cdot (9c_2^2\chi_2^2 +{\cred c_4}) \|v_0\|_{L^1(\Omega)}^4$, from (\ref{31.6}) we thus
  infer that
  \bas
    \frac{d}{dt} \io v^4 + \io v^4
    \le \frac{k_2(M) \ovz^4}{\tau^{10}} + c_5 e^{-\alpha_1 t}
    \qquad \mbox{for all } t\in (0,\tm).
  \eas
  Through Lemma \ref{lem21}, applicable here since $\alpha_1<1$, this entails that
  \bas
    \io v^4 \le \frac{k_2(M) \ovz^4}{\tau^{10}}
    + \bigg\{ \io v_0^4 + \frac{c_5}{1-\alpha_1} \bigg\} \cdot e^{-\alpha_1 t}
    \qquad \mbox{for all } t\in (0,\tm)
  \eas
  and hence completes the proof.
\qed
Yet concentrating on $v$ itself, we next resort to a semigroup-based argument once more to turn the above into an esimate
involving the norm in $L^\infty(\Omega)$.
\begin{lem}\label{lem32}
  One can find $\alpha>0$ in such a manner that for each $M>0$ there exists $K(M)>0$ such that
  if (\ref{init}) and (\ref{M}) are satisfied, then
  \be{32.1}
    \|v(\cdot,t)\|_{L^\infty(\Omega)}
    \le \frac{K(M) \ovz}{\tau^\frac{7}{2}} + C e^{-\alpha t}
    \qquad \mbox{for all } t\in (0,\tm)
  \ee
  with some $C=C(u_0,v_0,w_0)>0$.
\end{lem}
\proof
  A verification of this can be achieved in a way quite similar to that in Lemma \ref{lem24}:
  By the H\"older inequality as well as Lemma \ref{lem31} and Lemma \ref{lem22}, we see that with some $\alpha\in (0,1)$,
  given any $M>0$ we can find $k_1(M)>0$ such that if (\ref{init}) and (\ref{M}) hold, there exists $c_1=c_1(u_0,v_0,w_0)>0$
  fulfilling
  \bas
    \|v u_x\|_{L^\frac{4}{3}(\Omega)}
    \le \|v\|_{L^4(\Omega)} \|u_x\|_{L^2(\Omega)}
    \le \frac{k_1(M) \ovz}{\tau^\frac{7}{2}} + c_1 e^{-\alpha t}
    \qquad \mbox{for all } t\in (0,\tm).
  \eas
  Henceforth assuming (\ref{init}) and (\ref{M}), we combine this with known regularization features of the Neumann heat
  semigroup and (\ref{mass}) to see that with some positive constants $c_2$ and $c_3$ independent of $u_0, v_0$ and $w_0$ we
  have
  \bas
    \|v(\cdot,t)\|_{L^\infty(\Omega)}
    &=& \Bigg\| e^{t(\Delta-1)} v_0
    - \chi_2 \int_0^t e^{(t-s)(\Delta-1)} \partial_x \Big( v(\cdot,s) u_x(\cdot,s)\Big) ds
    + \int_0^t e^{(t-s)(\Delta-1)} v(\cdot,s) ds \Bigg\|_{L^\infty(\Omega)} \\
    &\le& e^{-t} \|e^{t\Delta} v_0\|_{L^\infty(\Omega)}
    + c_2 \int_0^t e^{-(t-s)} \Big(1+(t-s)^{-\frac{7}{8}}\Big) \|v(\cdot,s) u_x(\cdot,s)\|_{L^\frac{4}{3}(\Omega)} ds \\
    & & + c_2 \int_0^t e^{-(t-s)} \Big(1+(t-s)^{-\frac{1}{2}}\Big) \|v(\cdot,s)\|_{L^1(\Omega)} ds \\
    &\le& e^{-t} \|v_0\|_{L^\infty(\Omega)}
    + \frac{c_2 k_1(M) \ovz}{\tau^\frac{7}{2}} \int_0^t e^{-(t-s)} \Big(1+(t-s)^{-\frac{7}{8}}\Big) ds \\
    & & + c_1 c_2 \int_0^t e^{-(t-s)} \Big(1+(t-s)^{-\frac{7}{8}}\Big) e^{-\alpha s} ds \\
    & & + c_2 |\Omega| \ovz \int_0^t e^{-(t-s)} \Big(1+(t-s)^{-\frac{1}{2}}\Big) ds \\
    &\le& e^{-t} \|v_0\|_{L^\infty(\Omega)}
    + \frac{c_2 k_1(M) \ovz}{\tau^\frac{7}{2}} \int_0^\infty e^{-\sigma} (1+\sigma^{-\frac{7}{8}}) d\sigma \\
    & & + c_1 c_2 e^{-\alpha t} \int_0^\infty e^{-(1-\alpha)\sigma} (1+\sigma^{-\frac{7}{8}}) d\sigma \\
    & & + c_2 |\Omega| \ovz \int_0^\infty e^{-\sigma}(1+\sigma^{-\frac{1}{2}}) d\sigma
    \qquad \mbox{for all } t\in (0,\tm),
  \eas
  which readily yields (\ref{32.1}) due to the inequalities $\alpha<1$ and $\tau\le 1$.
\qed
As we are now in a position quite identical to that encountered immediately before Lemma \ref{lem22}, we can repeat
the argument thereof to finally derive the following gradient estimate for the crucial second solution component.
\begin{lem}\label{lem322}
  There exists $\alpha>0$ such that for all $M>0$ one can choose $K(M)>0$
  with the property that if (\ref{init}) and (\ref{M}) hold, the
  with some $C=C(u_0,v_0,w_0)>0$, we have
  \be{322.1}
    \io v_x^2(\cdot,t)
    \le \frac{K(M)}{\tau^{10}} + C e^{-\alpha t}
    \qquad \mbox{for all } t\in (0,\tm).
  \ee
\end{lem}
\proof
  The claimed inequality can be derived by means of an essentially verbatim copy of the argument from Lemma \ref{lem22},
  instead of referring to Lemma \ref{lem24} and Lemma \ref{lem29} now relying on Lemma \ref{lem32} and (\ref{22.2});
  we may therefore refrain from giving details here.
\qed
\mysection{Global existence and boundedness. Proof of Theorem \ref{theo34}}
Now asserting global extensibility of our solution actually reduces to a mere collection of our previously
obtained estimates, where at this stage neither any knowledge on the precise dependence thereof on $M$ or on $\ouz$ and
$\ovz$ is needed, nor do we rely on the exponentially decaying contributions to the above inequalities.
\begin{lem}\label{lem33}
  For all $u_0, v_0$ and $w_0$ fulfilling (\ref{init}), we have $\tm=\infty$, and furthermore we can find
  $C=C(u_0,v_0,w_0)>0$ such that (\ref{34.1}) holds.
\end{lem}
\proof
  In view of the extensibility criterion (\ref{ext}) from Lemma \ref{lem_loc}, for any fixed $(u_0,v_0,w_0)$ satisfying
  (\ref{init}) we may apply Lemma \ref{lem22}, Lemma \ref{lem322} and Lemma \ref{lem25} to $M:=\io u_0+\io v_0$ and $q:=2$ and
  thereby readily obtain that indeed $\tm$ cannot be finite, and that hence moreover (\ref{34.1}) is a consequence of
  (\ref{22.1}), (\ref{322.1}) and (\ref{25.1}).
\qed
In other words, we thereby already have derived our main result on global classical solvability in (\ref{0}):\abs
\proofc of Theorem \ref{theo34}. \quad
  We only need to combine Lemma \ref{lem33} with Lemma \ref{lem_loc}.
\qed
\mysection{Convergence for small values of $\min\{ \ouz, \ovz\}$. Proof of Theorem \ref{theo40}}
Next, in contrast to our development of the above existence statement,
our investigation of the large time behavior in (\ref{0}), as forming the objective of this section,
will considerably benefit from the more detailed information provided by our estimates from Section \ref{sect3}.
\subsection{Identifying a conditional energy functional}
The following lemma basically only
collects the essence of what will be needed from Section \ref{sect3} for our subsequent qualitative analysis.
\begin{lem}\label{lem35}
  Let $M>0$. Then there exists $K(M)>0$ with the property that if (\ref{init}) and (\ref{M}) hold, then one can
  find $t_0=t_0(u_0,v_0,w_0) \ge 0$ such that
  \be{35.1}
    \|u(\cdot,t)\|_{L^\infty(\Omega)} \le K(M) \ouz
    \qquad \mbox{for all } t>t_0
  \ee
  and
  \be{35.2}
    \|v(\cdot,t)\|_{L^\infty(\Omega)} \le K(M) \ovz
    \qquad \mbox{for all } t>t_0
  \ee
  as well as
  \be{35.3}
    \|w(\cdot,t)\|_{L^\infty(\Omega)} \le K(M)
    \qquad \mbox{for all } t>t_0.
  \ee
\end{lem}
\proof
  From Lemma \ref{lem24}, Lemma \ref{lem32} and Lemma \ref{lem26} we infer the existence of $\alpha>0$ such that whenever
  $M>0$, one can find $k_1(M)>0$ such that if (\ref{init}) and (\ref{M}) hold, we have
  \be{35.4}
    \|u(\cdot,t)\|_{L^\infty(\Omega)}
    \le k_1(M) \ouz + e^{-\alpha t}
    \qquad \mbox{for all } t>0
  \ee
  and
  \be{35.5}
    \|v(\cdot,t)\|_{L^\infty(\Omega)}
    \le k_1(M) \ovz + e^{-\alpha t}
    \qquad \mbox{for all } t>0
  \ee
  as well as
  \be{35.6}
    \|w(\cdot,t)\|_{L^\infty(\Omega)}
    \le k_1(M) + e^{-\alpha t}
    \qquad \mbox{for all } t>0,
  \ee
  where we note that in light of the fact that $\tm=\infty$ we now know that the number $\tau$ in Lemma \ref{lem32}
  actually satisfies $\tau=1$.
  For $t_0:=\frac{1}{\alpha} \cdot \ln_+ \frac{1}{k_1(M) \cdot \min\{\ouz,\ovz,1\}}$ and with $K(M):= 2k_1(M)$,
  the claimed inequalities now directly result from (\ref{35.4})-(\ref{35.6}).
\qed
Now a key toward our proof of stabilization can be found in the following observation on a genuine energy-type structure
in (\ref{0}) {\em when restricted to trajectories corresponding to initial data compatible with (\ref{M}) and (\ref{del1}).}
The presence of such conditional energy functionals, 
interpretable as a rigorous mathematical manifestation
of superlinear dependence on the unknown in the crucial nonlinearities, has been used in several studies
on asymptotic behavior in related chemotaxis problems in the recent few years (see e.g.~\cite{wang_xiang_yu},
\cite{chae_kang_lee_CPDE}, \cite{junping_shi_JDE},
\cite{win_ARMA}, \cite{qingshan_zhang_MANA2016} or also \cite{win_ct_sing_abs_eventual}     
for an incomplete collection); in comparison to most of these, the
seemingly most unique feature of the present situation consists in
that here it is possible to relax the smallness condition appearing
therein in such a substantial manner that in its remaining part it
merely reduces to a smallness assumption essentially equivalent to
(\ref{M})-(\ref{del1}):
\begin{lem}\label{lem36}
  Let $M>0$. Then there exists $\delta(M)>0$ such that if $u_0, v_0$ and $w_0$ are such that if beyond
  (\ref{init}) and (\ref{M}) we have
  \be{del}
    \bigg\{ \io u_0 \bigg\} \cdot \bigg\{ \io v_0\bigg\}^2 \le \delta(M),
  \ee
  then it is possible to find
  $b=b(u_0,v_0,w_0)>0$ and $t_0=t_0(u_0,v_0,w_0)>0$ with the property that
  \be{F}
    \F(t):=\io u(\cdot,t) \ln \frac{u(\cdot,t)}{\ouz}
    + b \io v(\cdot,t) \ln \frac{v(\cdot,t)}{\ovz}
    + \frac{\chi_1}{2\lambda} \io \frac{w_x^2(\cdot,t)}{w(\cdot,t)},
    \qquad t>0,
  \ee
  and
  \be{D}
    \D(t):=\frac{1}{2} \io \frac{u_x^2(\cdot,t)}{u(\cdot,t)}
    +\frac{b}{2} \io \frac{v_x^2(\cdot,t)}{v(\cdot,t)}
    + \frac{\chi_1 \mu}{4\lambda} \io \frac{w_x^2(\cdot,t)}{w(\cdot,t)},
    \qquad t>0,
  \ee
  satisfy
  \be{36.2}
    \F'(t) \le -\D(t)
    \qquad \mbox{for all } t>t_0.
  \ee
\end{lem}
\proof
  Given $M>0$, we first apply Lemma \ref{lem35} to fix $k_1(M)>0$ with the property that for any choice of $(u_0,v_0,w_0)$
  complying with (\ref{init}) and (\ref{M}) we can find $t_0(u_0,v_0,w_0)\ge 0$ such that for all $t>t_0$,
  \be{36.3}
    \|u(\cdot,t)\|_{L^\infty(\Omega)}
    \le k_1(M) \ouz,
    \quad
    \|v(\cdot,t)\|_{L^\infty(\Omega)}
    \le k_1(M) \ovz
    \quad \mbox{and} \quad
    \|w(\cdot,t)\|_{L^\infty(\Omega)}
    \le k_1(M),
  \ee
  and we thereupon claim that the intended conclusion holds if we let
  \be{36.4}
    \delta(M):=\frac{\mu |\Omega|^3}{8\chi_1 \chi_2^2 \lambda k_1^4(M)}.
  \ee
  To see this, we fix any $(u_0,v_0,w_0)$ fulfilling (\ref{init}) and (\ref{M}) as well as (\ref{del}), and abbreviating
  $t_0:=t_0(u_0,v_0,w_0)$, $L_u:=k_1(M)\ouz, L_v:=k_1(M)\ovz$ and $L_w:=k_1(M)$ we infer from (\ref{36.4}) that it is possible
  to pick $b=b(u_0,v_0,w_0)>0$ such that
  \be{36.44}
    \frac{4\chi_1 \lambda  L_v L_w}{\mu} \le b \le \frac{1}{2\chi_2^2 L_u L_v}.
  \ee
  We then let $\F$ and $\D$ be as accordingly defined through (\ref{F}) and (\ref{D}), and in order to verify (\ref{36.2})
  we integrate by parts in (\ref{0}) and use (\ref{mass}) to compute
  \be{36.5}
    \frac{d}{dt} \io u\ln \frac{u}{\ouz}
    = \frac{d}{dt} \io u\ln u
    = - \io \frac{u_x^2}{u} + \chi_1 \io u_x w_x
  \ee
  and
  \be{36.6}
    \frac{d}{dt} \io v\ln \frac{v}{\ovz}
    = \frac{d}{dt} \io v\ln v
    = - \io \frac{v_x^2}{v}
    + \chi_2 \io u_x v_x
  \ee
  as well as
  \bea{36.7}
    \frac{d}{dt} \io \frac{w_x^2}{w}
    &=& 2 \io \frac{w_x}{w} \cdot \Big\{ dw_{xxx} - \lambda u_x w - \lambda uw_x - \lambda v_x w - \lambda vw_x
        - \mu w_x \Big\} \nn\\
    & & - \io \frac{w_x^2}{w^2} \cdot \Big\{ dw_{xx} -\lambda uw -\lambda vw - \mu w + r \Big\} \nn\\
    &=& -2d \io \frac{w_{xx}^2}{w}
    + d\io \frac{w_x^2 w_{xx}}{w} \nn\\
    & & -2\lambda \io u_x w_x
    - 2\lambda \io v_x w_x
    - \lambda \io \frac{u}{w} w_x^2
    - \lambda \io \frac{v}{w} w_x^2
    - \mu \io \frac{w_x^2}{w}
  \eea
  for $t>0$.
  Here once more integrating by parts we see that
  \bas
    \io \frac{w_x^4}{w^3}
    = -\frac{1}{2} \io \Big( \frac{1}{w^2}\Big)_x w_x^3
    = \frac{3}{2} \io \frac{w_x^2 w_{xx}}{w^2}
    \qquad \mbox{for all } t>0,
  \eas
  which by the Cauchy-Schwarz inequality firstly entails that
  \bas
    \io \frac{w_x^4}{w^3} \le \frac{3}{2} \bigg\{ \io \frac{w_{xx}^2}{w} \bigg\}^\frac{1}{2} \cdot
    \bigg\{ \io \frac{w_x^4}{w^3}\bigg\}^\frac{1}{2}
    \qquad \mbox{for all } t>0
  \eas
  and hence
  \bas
    \io \frac{w_x^4}{w^3} \le \frac{9}{4} \io \frac{w_{xx}^2}{w}
    \qquad \mbox{for all } t>0,
  \eas
  and which, as a consequence, secondly shows that thus in (\ref{36.7}) we can estimate
  \bas
    -2d \io \frac{w_{xx}^2}{w} + d\io \frac{w_x^2 w_{xx}}{w}
    = -2d\io \frac{w_{xx}^2}{w}
    + \frac{2d}{3} \io \frac{w_x^4}{w^3}
    \le -\frac{d}{2} \io \frac{w_{xx}^2}{w} \le 0
    \qquad \mbox{for all } t>0.
  \eas
  Upon combining (\ref{36.5})-(\ref{36.7}) and neglecting two further well-signed summands, we therefore obtain that
  \be{36.8}
    \F'(t) + \io \frac{u_x^2}{u} + b \io \frac{v_x^2}{v} + \frac{\chi_1 \mu}{2\lambda} \io \frac{w_x^2}{w}
    \le b \chi_2 \io u_x v_x - \chi_1 \io v_x w_x
    \qquad \mbox{for all } t>0.
  \ee
  Here by Young's inequality and (\ref{36.3}),
  \bea{36.9}
    b\chi_2 \io u_x v_x
    &\le& \frac{1}{2} \io \frac{u_x^2}{u}
    + \frac{b^2 \chi_2^2}{2} \io u v_x^2 \nn\\
    &\le& \frac{1}{2} \io \frac{u_x^2}{u}
    + \frac{b^2 \chi_2^2}{2} \|u\|_{L^\infty(\Omega)} \|v\|_{L^\infty(\Omega)} \io \frac{v_x^2}{v} \nn\\
    &\le& \frac{1}{2} \io \frac{u_x^2}{u}
    + \frac{b^2 \chi_2^2 L_u L_v}{2} \io \frac{v_x^2}{v} \nn\\
    &\le& \frac{1}{2} \io \frac{u_x^2}{u}
    + \frac{b}{4} \io \frac{v_x^2}{v}
    \qquad \mbox{for all } t>t_0,
  \eea
  because thanks to the right inequality in (\ref{36.44}) we know that
  \bas
    \frac{\frac{b^2 \chi_2^2 L_u L_v}{2}}{\frac{b}{4}}
    = 2b\chi_2^2 L_u L_v \le 1.
  \eas
  Likewise, Young's inequality together with (\ref{36.3}) moreover shows that
  \bea{36.10}
    -\chi_1 \io v_x w_x
    &\le& \frac{b}{4} \io \frac{v_x^2}{v}
    + \frac{\chi_1^2}{b} \io vw_x^2 \nn\\
    &\le& \frac{b}{4} \io \frac{v_x^2}{v}
    + \frac{\chi_1^2}{b} \|v\|_{L^\infty(\Omega)} \|w\|_{L^\infty(\Omega)} \io \frac{w_x^2}{w} \nn\\
    &\le& \frac{b}{4} \io \frac{v_x^2}{v}
    + \frac{\chi_1^2 L_v L_w}{b} \io \frac{w_x^2}{w} \nn\\
    &\le& \frac{b}{4} \io \frac{v_x^2}{v}
    + \frac{\chi_1 \mu}{4\lambda} \io \frac{w_x^2}{w}
    \qquad \mbox{for all } t>t_0,
  \eea
  since by the left restriction in (\ref{36.44}),
  \bas
    \frac{\frac{\chi_1^2 L_v L_w}{b}}{\frac{\chi_1 \mu}{4\lambda}}
    = \frac{4\chi_1 \lambda L_v L_w}{b\mu} \le 1.
  \eas
  It thus remains to insert (\ref{36.9}) and (\ref{36.10}) into (\ref{36.8}) to end up with (\ref{36.2}).
\qed
\subsection{Exponential convergence. Proof of Theorem \ref{theo40}}
A first and rather immediate consequence of (\ref{36.2}) when combined with well-known inequalities of logarithmic Sobolev and
Csisz\'ar-Kullback type yields convergence already at algebraic rates, but yet with respect to spatial $L^1$ norms only.
\begin{lem}\label{lem37}
  Given $M>0$, let $\delta(M)>0$ be as in Lemma \ref{lem36}, and suppose that (\ref{init}), (\ref{M}) and (\ref{del}) hold.
  Then there exist $C=C(u_0,v_0,w_0)>0$ and $\alpha=\alpha(u_0,v_0,w_0)>0$ such that
  \be{37.1}
    \|u(\cdot,t)-\ouz\|_{L^1(\Omega)}
    + \|v(\cdot,t)-\ovz\|_{L^1(\Omega)}
    \le C e^{-\alpha t}
    \qquad \mbox{for all } t>0.
  \ee
\end{lem}
\proof
  Given $(u_0,v_0,w_0)$ such that (\ref{init}), (\ref{M}) and (\ref{del}) hold, we take $b=b(u_0,v_0,w_0)>0$ and
  $t_0=t_0(u_0,v_0,w_0)>0$ as provided by Lemma \ref{lem36}, and recall that according to a logarithmic Sobolev inequality
  (\cite{gross}, \cite{rothaus})
  and (\ref{mass}) there exists $c_1>0$ such that
  \bas
    \io u\ln \frac{u}{\ouz} +  b\io v\ln \frac{v}{\ovz}
    \le c_1 \cdot \bigg\{ \frac{1}{2} \io \frac{u_x^2}{u} + \frac{b}{2} \io \frac{v_x^2}{v}\bigg\}
    \qquad \mbox{for all } t>0.
  \eas
  Writing $c_2:=\max\{c_1,\frac{2}{\mu}\}$, for $\F$ and $\D$ as in (\ref{F}) and (\ref{D}) we thus obtain that
  \bas
    \F(t) \le c_2 \D(t)
    \qquad \mbox{for all } t>0,
  \eas
  so that (\ref{36.2}) implies the autonomous ODI
  \bas
    \F'(t) \le -\frac{1}{c_2} \F(t)
    \qquad \mbox{for all } t>t_0.
  \eas
  Upon integration, this entails that
  \bas
    \F(t) \le \F(t_0) e^{-\frac{t-t_0}{c_2}}
    \qquad \mbox{for all } t>t_0
  \eas
  and thereby establishes (\ref{37.1}) with suitably large $C>0$ and $\alpha:=\frac{1}{2c_2}$, because
  according to a Csisz\'ar-Kullback inequality (\cite{csiszar}, \cite{arnold}) and (\ref{mass}) we can find $c_3>0$ fulfilling
  \bas
    \|u(\cdot,t)-\ouz\|_{L^1(\Omega)}^2
    \le c_3\io u(\cdot,t) \ln \frac{u(\cdot,t)}{\ouz}
    \le c_3 \F(t)
  \eas
  and
  \bas
    \|v(\cdot,t)-\ovz\|_{L^1(\Omega)}^2
    \le c_3\io v(\cdot,t) \ln \frac{v(\cdot,t)}{\ovz}
    \le \frac{c_3}{b} \F(t)
  \eas
  for all $t>0$.
\qed
Thanks to the temporally uniform $H^1$ bounds for both $u$ and $v$ known from Theorem \ref{theo34},
a straightforward interpolation finally asserts exponential convergence also with respect to $L^\infty$ norms.
\begin{lem}\label{lem38}
  Under the assumptions of Lemma \ref{lem37}, one can find
  $C=C(u_0,v_0,w_0)>0$ and $\alpha=\alpha(u_0,v_0,w_0)>0$ such that
  \be{38.1}
    \|u(\cdot,t)-\ouz\|_{L^\infty(\Omega)}
    \le C e^{-\alpha t}
    \qquad \mbox{for all } t>0
  \ee
  and
  \be{38.2}
    \|v(\cdot,t)-\ovz\|_{L^\infty(\Omega)}
    \le C e^{-\alpha t}
    \qquad \mbox{for all } t>0.
  \ee
\end{lem}
\proof
  By means of a Gagliardo-Nirenberg interpolation, we can find $c_1>0$ such that
  \bas
    \|\varphi-\overline{\varphi}\|_{L^\infty(\Omega)}
    \le c_1\|\varphi_x\|_{L^2(\Omega)}^\frac{2}{3} \|\varphi-\overline{\varphi}\|_{L^1(\Omega)}^\frac{1}{3}
    \qquad \mbox{for all } \varphi\in W^{1,2}(\Omega),
  \eas
  whereas Theorem \ref{theo34} warrants the existence of $c_2>0$ such that
  \bas
    \|u_x(\cdot,t)\|_{L^2(\Omega)} \le c_2
    \quad \mbox{and} \quad
    \|v_x(\cdot,t)\|_{L^2(\Omega)} \le c_2
    \qquad \mbox{for all } t>0.
  \eas
  In view of (\ref{mass}), we therefore obtain that for all $t>0$,
  \bas
    \|u(\cdot,t)-\ouz\|_{L^\infty(\Omega)}
    + \|v(\cdot,t)-\ouz\|_{L^\infty(\Omega)}
    \le c_1 c_2^\frac{2}{3} \|u(\cdot,t)-\ouz\|_{L^1(\Omega)}^\frac{1}{3}
    + c_1 c_2^\frac{2}{3} \|v(\cdot,t)-\ovz\|_{L^1(\Omega)}^\frac{1}{3},
  \eas
  which due to Lemma \ref{lem37} implies both (\ref{38.1}) and (\ref{38.2}).
\qed
On the basis of a final testing procedure independent from the above, again combined with an interpolation argument
of the above flavor, we can lastly derive exponential and spatially uniform stabilization also in the third solution component.
\begin{lem}\label{lem39}
  Suppose that the hypotheses of Lemma \ref{lem37} are satisfied. Then there exist
  $C=C(u_0,v_0,w_0)>0$ and $\alpha=\alpha(u_0,v_0,w_0)>0$ such that with $\ws\ge 0$ given by (\ref{ws}) we have
  \be{39.1}
    \|w(\cdot,t)-\ws\|_{L^\infty(\Omega)}
    \le C e^{-\alpha t}
    \qquad \mbox{for all } t>0.
  \ee
\end{lem}
\proof
  Using that by definition of $\ws$ we have $\{ \lambda(\ouz+\ovz) + \mu\} \ws=r$, upon testing the third equation in (\ref{0})
  by $w-\ws$ we obtain that for all $t>0$,
  \bea{39.2}
    \frac{1}{2} \frac{d}{dt} \io (w-\ws)^2
    + d\io w_x^2
    &=& \io \Big\{ -\lambda (u+v)w - \mu w + r\Big\} \cdot (w-\ws) \nn\\
    &=& \io \Big\{ -\lambda (\ouz+\ovz) w - \mu w + r\Big\} \cdot (w-\ws) \nn\\
    & & - \lambda  \io (u-\ouz) w(w-\ws)
    - \lambda \io (v-\ovz) w(w-\ws) \nn\\
    &=& - c_1 \io (w-\ws)^2 \nn\\
    & & - \lambda  \io (u-\ouz) w(w-\ws)
    - \lambda \io (v-\ovz) w(w-\ws)
  \eea
  with $c_1:=\lambda (\ouz+\ovz) +\mu>0$.
  Here since $\|w(\cdot,t)\|_{L^\infty(\Omega)} \le c_2:=\frac{r}{\mu} + \|w_0\|_{L^\infty(\Omega)}$ for all $t>0$ by
  Lemma \ref{lem26}, Young's inequality shows that
  \bas
    -\lambda \io (u-\ouz) w(w-\ws)
    &\le& \frac{c_1}{4} \io (w-\ws)^2
    + \frac{\lambda^2}{c_1} \io (u-\ouz)^2 w^2 \nn\\
    &\le& \frac{c_1}{4} \io (w-\ws)^2
    + \frac{c_2^2 \lambda^2}{c_1} \io (u-\ouz)^2
    \qquad \mbox{for all } t>0,
  \eas
  and similarly estimating the rightmost summand in (\ref{39.2}) we altogether infer that
  \be{39.3}
    \frac{d}{dt} \io (w-\ws)^2
    + c_1 \io (w-\ws)^2
    \le c_3 \io (u-\ouz)^2
    + c_3 \io (v-\ovz)^2
    \qquad \mbox{for all } t>0
  \ee
  if we let $c_3:=\frac{2c_2^2 \lambda^2}{c_1}$.
  As Lemma \ref{lem38} provides $c_4>0$ and $\alpha_1\in (0,c_1)$ fulfilling
  \bas
    \io \Big(u(\cdot,t)-\ouz\Big)^2
    + \io \Big(v(\cdot,t)-\ovz\Big)^2
    \le c_4 e^{-\alpha_1 t}
    \qquad \mbox{for all } t>0,
  \eas
  through e.g.~Lemma \ref{lem21} we readily conclude from (\ref{39.3}) that
  \bas
    \io \Big(w(\cdot,t)-\ws\Big)^2 \le c_5 e^{-\alpha_1 t}
    \qquad \mbox{for all } t>0
  \eas
  with $c_5:=\io (w_0+\ws)^2 + \frac{c_4}{c_1-\alpha_1}$.\abs
  Now since Lemma \ref{lem33} in conjunction with Lemma \ref{lem25} asserts the existence of $c_6>0$ such that
  $\|w(\cdot,t)-\ws\|_{W^{1,2}(\Omega)} \le c_6$ for all $t>0$, and hence by the Gagliardo-Nirenberg inequality we can find
  $c_7>0$ such that
  \bas
    \|w(\cdot,t)-\ws\|_{L^\infty(\Omega)}
    \le c_7\|w(\cdot,t)-\ws\|_{W^{1,2}(\Omega)}^\frac{1}{2} \|w(\cdot,t)-\ws\|_{L^2(\Omega)}^\frac{1}{2}
    \qquad \mbox{for all } t>0,
  \eas
  this implied that
  \bas
    \|w(\cdot,t)-\ws\|_{L^\infty(\Omega)}
    \le c_5^\frac{1}{4} c_6^\frac{1}{2} c_7 e^{-\frac{\alpha_1}{4} t}
    \qquad \mbox{for all } t>0
  \eas
  and hence entails (\ref{39.1}).
\qed
Our main statements on temporal asymptotics thereby become almost evident:\abs
\proofc of Theorem \ref{theo40}. \quad
  For $M>0$ we take $\delta(M)>0$ as provided by Lemma \ref{lem36} and let $\eps(M)>0$ be small enough such that
  $\eps(M) \le \frac{\delta(M)}{M^2}$ and $\eps(M) \le \sqrt{\frac{\delta}{M)}{M}}$.
  It can then readily be verified that assuming (\ref{M}) together with (\ref{del1}) entails (\ref{del}), so that
  for completing the proof it is sufficient to collect the statements from
  Lemma \ref{lem38} and Lemma \ref{lem39}.
\qed
\mysection{Appendix: Two statements on ODE comparison}
Let us finally state two elementary results of quite straightforward ODE comparison arguments, in view of our above applications
with particular focus on the respective dependence on the parameters appearing therein.
We begin with a simple observation that has been used in the proofs of Lemma \ref{lem23}, Lemma \ref{lem31}
and Lemma \ref{lem39}.
\begin{lem}\label{lem21}
  Let $\kappa>0, a\ge 0, b\ge 0$ and $\alpha>0$ be such that $\alpha<\kappa$, and suppose that
  $y\in C^0([0,T)) \cap C^1((0,T))$ is a nonnegative function satisfying
  \be{21.1}
    y'(t) + \kappa y(t) \le a  e^{-\alpha t} + b
    \qquad \mbox{for all } t\in (0,T)
  \ee
  with some $T\in (0,\infty]$. Then
  \be{21.2}
    y(t) \le \Big( y(0) + \frac{a}{\kappa-\alpha}\Big) e^{-\alpha t} + \frac{b}{\kappa}
    \qquad \mbox{for all } t\in (0,T).
  \ee
\end{lem}
\proof
  We let $c_1:=y(0)+\frac{a}{\kappa-\alpha}$ and $\oy(t):=c_1 e^{-\alpha t} + \frac{b}{\kappa}$ for $t\ge 0$.
  Then clearly $\oy(0)>y(0)$, and since moreover
  $\oy'(t) + \kappa \oy(t) - a e^{-\alpha t} - b=\{(\kappa-\alpha)c_1 -a\} e^{-\alpha t} =(\kappa-\alpha)y(0)e^{-\alpha t}\ge 0$ for all $t>0$
  thanks to the nonnegativity of $y$ and our assumption that $\alpha<\kappa$,
  the inequality in (\ref{21.2}) results from a comparison argument.
\qed
Our second statement in this direction, as used in crucial places in Lemma \ref{lem29} and Lemma \ref{lem22},
merely imposes some hypothesis on temporal averages of the respective force term, and therefore requires a slightly
more subtle argument:
\begin{lem}\label{lem211}
  Let $\kappa>0, a\ge 0, b\ge 0$ and $\alpha\in (0,\kappa)$, and assume that with some $T\in (0,\infty]$ and $\tau\in (0,T)$
  such that $\tau\le 1$, the nonnegative functions
  $y\in C^0([0,T)) \cap C^1((0,T))$ and $f\in L^1_{loc}([0,T))$ are such that
  \be{211.1}
    y'(t) + \kappa y(t) \le f(t)
    \qquad \mbox{for all } t\in (0,T)
  \ee
  and
  \be{211.2}
    \int_t^{t+\tau} f(s) ds \le ae^{-\alpha t} + b
    \qquad \mbox{for all } t \in [0,T-\tau).
  \ee
  Then
  \be{211.3}
    y(t) \le \bigg\{ \Big( y(0)+  a+ \frac{a}{\kappa-\alpha} + b\Big) \cdot \frac{e^\alpha}{\tau}
        + ae^\alpha \bigg\} \cdot e^{-\alpha t}
    + \frac{b}{\kappa\tau} + b
    \qquad \mbox{for all } t\in (0,T).
  \ee
\end{lem}
\proof
  Letting $z(t):=\int_t^{t+\tau} y(s) ds$ for $t\in [0, T-\tau)$, by two integrations of (\ref{211.1}),
  from (\ref{211.2}) we obtain that firstly
  \be{211.4}
    z'(t) + \kappa z(t) \le  ae^{-\alpha t} + b
    \qquad \mbox{for all } t\in (0,T-\tau),
  \ee
  and that secondly,
  \be{211.6}
    y(t) \le y(t_0) + \int_{t_0}^t f(s) ds
    \le y(t_0) +  a e^{-\alpha t_0} + b
    \qquad \mbox{for all $t_0\in [0,T)$ and $t\in (t_0,T)$ such that } t\le t_0+\tau,
  \ee
  whence in particular
  \be{211.66}
    y(t) \le y(0) +  a + b
    \qquad \mbox{for all } t\in  (0,\tau]
  \ee
  As thus
  \bas
    z(0) = \int_0^\tau y(s) ds \le y(0) + a+b
  \eas
  due to our assumption that $\tau\le 1$,
  by applying Lemma \ref{lem21} to (\ref{211.4}) we see that
  \bas
    z(t)
    &\le& \Big( z(0) + \frac{a}{\kappa-\alpha} \Big) e^{-\alpha t} + \frac{b}{\kappa} \nn\\
    &\le& \Big(y(0)+ a+\frac{a}{\kappa-\alpha}+b\Big) e^{-\alpha t} + \frac{b}{\kappa}
    \qquad \mbox{for all } t\in (0,T-\tau),
  \eas
  which implies that
  \be{211.7}
    z(t-\tau) \le \Big(y(0)+a+\frac{a}{\kappa-\alpha}+b\Big) e^{-\alpha (t-\tau)}  + \frac{b}{\kappa}
        \qquad \mbox{for all } t\in (\tau, T).
  \ee
  In view of the definition of $z(t-\tau)$, this especially means that whenever  $t\in
  (\tau,T)$, we can find $t_0(t) \in (t-\tau,t)$ fulfilling
  \bas
    y(t_0(t))
    \le \frac{1}{\tau} \cdot \Big( y(0)+ a + \frac{a}{\kappa-\alpha} +b\Big) \cdot e^{-\alpha (t-\tau)}
    +\frac{b}{\kappa  \tau}.
  \eas
  Again employing (\ref{211.6}), and moreover using the that $e^{-\alpha t_0}<e^{-\alpha(t-\tau)}$ due to the inclusion
  $t_0\in (t-\tau, t)$, we conclude that for any such $t$,
  \bas
    y(t)
    &\le& y(t_0(t)) +  a e^{-\alpha t_0} +b \\
    &\le& \bigg\{ \frac{1}{\tau} \Big( y(0)+ a + \frac{a}{\kappa-\alpha} +b\Big) e^{\alpha\tau}
        + a e^{\alpha\tau} \bigg\} \cdot e^{-\alpha t}
    + \frac{b}{\kappa \tau} + b,
  \eas
  from which (\ref{211.3}) immediately follows once more due to the inequality $\tau\le 1$.\abs
  If $t\in (0,\tau]$, however, we infer (\ref{211.3}) directly from (\ref{211.66}), because clearly
  \bas
    y(0)+ a + b
    \le \frac{1}{\tau} \Big(y(0)+  a+b\Big) e^{\alpha\tau} e^{-\alpha t}
  \eas
  for all $t\in (0,\tau]$.
\qed

\vspace*{10mm}
{\bf Acknowledgment.}
  Y. Tao acknowledges support of the {\em National Natural Science Foundation of China (No.
  11571070)}.
  M. Winkler was supported by the {\em Deutsche Forschungsgemeinschaft} within the project
  {\em Analysis of chemotactic cross-diffusion in complex
  frameworks}.
\end{document}